\documentclass[letterpaper]{article}

\usepackage{arxiv}

\usepackage[utf8]{inputenc}
\usepackage{float}
\usepackage{enumitem}
\usepackage[english]{babel}
\usepackage{amssymb}
\usepackage{amsmath,mathtools}
\usepackage{algorithmicx}
\usepackage{algorithm}
\usepackage{algpseudocode}
\usepackage{xcolor}
\usepackage{graphicx}
\usepackage{multirow}
\usepackage{overpic}
\usepackage{psfrag}
\usepackage{bm}
\usepackage{url}
\usepackage{epsfig}
\usepackage{marvosym}
\usepackage{verbatim}
\usepackage{subcaption}
\usepackage{caption}
\usepackage{hyperref}
\hypersetup{colorlinks=true,linkcolor=red,citecolor=blue}
\usepackage{ulem}
\usepackage{lmodern}
\usepackage{booktabs} 
\usepackage{moresize} 
\usepackage{siunitx}
\usepackage{tikz}
\usepackage{pgfplots}
\pgfplotsset{compat=1.14}
\usepgfplotslibrary{fillbetween}
\usetikzlibrary{shapes, arrows, arrows.meta, calc, positioning, decorations,
decorations.pathmorphing, decorations.text, spy}
\pgfplotsset{translate gnuplot=true}
\usepackage{etex} 								
\usepgfplotslibrary{external} 						
\tikzsetexternalprefix{ext/}							
\usepackage{lineno}

\definecolor{myblue1}	{RGB}{0,177,234}				
\definecolor{myblue2}	{RGB}{76,200,239}				
\definecolor{myblue3}	{RGB}{127,215,244}				
\definecolor{myblue4}	{RGB}{178,231,248}				
\definecolor{myblue5}	{RGB}{198,251,255}				
\definecolor{mybluegray1}{RGB}{0,127,167}				
\definecolor{mybluegray2}{RGB}{76,165,193}				
\definecolor{mybluegray3}{RGB}{127,191,211}				
\definecolor{mybluegray4}{RGB}{178,216,228}				
\definecolor{mygray1}	{RGB}{76,84,93}				
\definecolor{mygray2}	{RGB}{129,135,141}				
\definecolor{mygray3}	{RGB}{165,169,174}				
\definecolor{mygray4}	{RGB}{201,203,206}				
\definecolor{myorange1}	{RGB}{255,126,46}				
\definecolor{myorange2}	{RGB}{255,164,108}				
\definecolor{myorange3}	{RGB}{255,190,150}				
\definecolor{myorange4}	{RGB}{255,216,192}				
\definecolor{mypurple1}	{RGB}{89,89,171}				
\definecolor{mypurple4}	{RGB}{189,189,231}				
\definecolor{mybrown1}	{RGB}{205,133,63}				
\definecolor{mybrown4}	{RGB}{205,183,123}				
\definecolor{mygreen1}	{RGB}{60,179,113}				
\definecolor{mygreen4}	{RGB}{144,238,144}				

\pgfplotsset{
    colormap={custom_map}{[5pt]
            rgb255(0pt)=(255,126,46);
            rgb255(500pt)=(255,190,150);
            rgb255(1000pt)=(0,177,234);
            rgb255(1500pt)=(127,215,244);
    },
}

\newcommand\ie{\textit{i.e.,}~}
\newcommand\eg{{e.g.,}~}

\renewcommand{\emph}[1]{\textit{#1}}
\setenumerate{label=\textcolor{mybluegray1}{$\bm{\diamond}$},itemsep=0.0pt,topsep=5pt}

\algnewcommand\algorithmicinitialize{\textbf{initialize}}
\algnewcommand\Initialize[1]{\State\algorithmicinitialize\ #1}

\algnewcommand\algorithmicsample{\textbf{sample}}
\algnewcommand\Sample[1]{\State\algorithmicsample\ #1}

\algnewcommand\algorithmicmap{\textbf{map}}
\algnewcommand\Map[1]{\State\algorithmicmap\ #1}

\algnewcommand\algorithmicprovide{\textbf{provide}}
\algnewcommand\Provide[1]{\State\algorithmicprovide\ #1}

\algnewcommand\algorithmicretrieve{\textbf{retrieve}}
\algnewcommand\Retrieve[1]{\State\algorithmicretrieve\ #1}

\algnewcommand\algorithmiccompute{\textbf{compute}}
\algnewcommand\Compute[1]{\State\algorithmiccompute\ #1}

\algnewcommand\algorithmicshuffle{\textbf{shuffle}}
\algnewcommand\Shuffle[1]{\State\algorithmicshuffle\ #1}

\algnewcommand\algorithmicgenerate{\textbf{generate}}
\algnewcommand\Generate[1]{\State\algorithmicgenerate\ #1}

\algnewcommand\algorithmicupdate{\textbf{update}}
\algnewcommand\Update[1]{\State\algorithmicupdate\ #1}

\title{Policy-based optimization: single-step policy gradient method seen as an evolution strategy}

\author{
	J. Viquerat\thanks{Corresponding author}\\
	MINES Paristech, CEMEF\\
	PSL - Research University\\
	\texttt{jonathan.viquerat@mines-paristech.fr}\\
\And
	R. Duvigneau\\
	INRIA Sophia Antipolis M\'editerran\'ee\\
	ACUMES project-team
\And
	P. Meliga\\
	MINES Paristech, CEMEF\\
	PSL - Research University
\And
	A. Kuhnle\\
	University of Cambridge
\And
	E. Hachem\\
	MINES Paristech, CEMEF\\
	PSL - Research University
}

\begin{document}
\newgeometry{left=3cm,right=3cm,top=3cm,bottom=2.5cm}
\maketitle

\begin{abstract} 
This research reports on the recent development of black-box optimization methods based on single-step deep reinforcement learning (DRL) and their conceptual similarity to evolution strategy (ES) techniques. It formally introduces 
policy-based optimization (PBO), a policy-gradient-based optimization algorithm that relies on a policy network to describe the density function of its forthcoming evaluations, and uses covariance estimation to steer the policy improvement process in the right direction. The specifics of the PBO algorithm are detailed, and the connection to evolutionary strategies is discussed. Relevance is assessed by benchmarking PBO against classical ES techniques on analytic functions minimization problems, and by optimizing various parametric control laws intended for the Lorenz attractor.
Given the scarce existing literature on the topic, this contribution definitely establishes PBO as a valid, versatile black-box optimization technique, and opens the way to multiple future improvements building on the inherent flexibility of the neural networks approach.
\end{abstract}


\section{Introduction}

During the past decade, machine learning methods, and more specifically deep neural network (DNN), have achieved great success in a wide variety of domains. State-of-the-art neural network architectures have reached astonishing performance levels in a variety of tasks, \ie image classification \cite{rawat2017,khan2020}, speech recognition \cite{nassif2019} or generative tasks \cite{gui2020}, to name a few. With generalized access to GPU computational resources through cheaper hardware or cloud computing, such advances have opened the path to a revolution of the reference methods in these domains, at both academic and industrial levels. 

With neural networks quickly becoming pervasive in a broad range of domains, significant progress has been made in solving challenging decision-making problems by deep reinforcement learning (DRL), an advanced branch of machine learning that couples DNNs and reinforcement learning (RL) algorithms. 
The ability to use high-dimensional state spaces and to exploit the feature extraction capabilities of DNNs has proven decisive to lift the obstacles that had long hindered classical RL methods. In return, this yielded unprecedented efficiency in games \cite{mnih2013,silver2017,openAI2018} and in several scientific disciplines such as robotics \cite{pinto2017}, language processing \cite{bahdanau2016}, although a tremendous potential also exists for applying DRL to real-life applications, including autonomous cars \cite{kendall2018,bewley2018} or data center cooling \cite{googleDataCenter2018}. 

Although neural networks are regularly used in optimization problems, they are most often exploited as trained surrogates for the actual objective function \cite{villarrubia2018,schweidtmann2019}, and have long been mostly left out of the central optimization process, with the exception of a handful of studies \cite{andrychowicz2016,yan2019,li2020}.
This trend has been changing lately, as tweaked versions of classical DRL policy gradient algorithms have began to be used as black-box optimizers \cite{viquerat2021}, the underlying idea being that a DRL agent can learn to map the same initial state to an optimal set by exploiting "single-step episodes", if the policy to be learnt is independent of state (hence, by extension, single-step DRL).
Such an approach has been speculated to hold a high potential for reliable optimization of complex systems. Nonetheless, it remains to be analyzed in full depth, as feasibility has just been assessed in 
a computational fluid mechanics (CFD) context, including shape optimization \cite{viquerat2021}, drag reduction \cite{ghraieb2020} and conjugate heat transfer control \cite{hachem2020} (a similar concept of "stateless DRL" has been early sketched in \cite{hmlinen2018} for validation purpose, but was not pursued).

This research formally introduces policy-based optimization (PBO), a novel single-step DRL algorithm that shares strong similarities with evolution strategies (ES). The objective is twofold: first, to deliver several major improvements over our previous PPO-1 algorithm, by adopting key heuristics from the covariance matrix adaptation evolution strategy (CMA-ES); second, to shape the capabilities of the method and to provide performance comparison against canonical ES algorithms on textbook and applied cases. 
An additional novelty lies in the generation of valid covariance matrices from neural network outputs, using the hypersphere decomposition method. To do so, three separate neural networks are exploited to learn the mean, variance and correlation parameters of a multivariate normal search distribution, yielding a powerful, flexible optimization method (without anticipating the results, PBO compares very well with CMA-ES, which is all the more promising since new algorithms cannot be expected to reach right away the level of performance of their more established counterparts). In comparison, PPO-1 updates the mean and variance (the same for all variables) from a single neural network, which can prematurely shrink the exploration variance. 
While there have been previous attempts to similarly improve the convergence properties of classical DRL algorithm by drawing inspiration from ES \cite{hmlinen2018}, our literature review did not reveal any other study considering the generation of valid full covariance matrices from neural network outputs.
 
The organization of the remaining of the paper is as follows: section \ref{section:pg-es} provides the needed background on policy gradient reinforcement learning and evolutionary strategies. 
The policy-based optimization (PBO) algorithm is introduced in section \ref{section:pbo}, 
where we thoroughly examine how to generate valid full covariance matrices from neural network outputs,
and point out key similarities and differences with respect to CMA-ES.\footnote{ The base code used to produce all results documented in this paper is available via a dedicated github repository \cite{pbo}. }
In section \ref{section:analytic}, PBO is benchmarked against standard ES algorithms on textbook optimization problems of 
analytic functions minimization. Finally, section \ref{section:lorenz} uses PBO to optimize a parametric control law for the Lorenz attractor, a well-known reduced version of Rayleigh--Bénard convection.

\section{Preliminaries}
\label{section:pg-es}

\subsection{Neural networks}

A neural network (NN) is a collection of artificial neurons, \ie connected computational units that can be trained to approximate arbitrarily well the mapping function between input and output spaces \cite{hornik1989}. Each connection provides the output of a neuron as an input to another neuron. Each neuron performs a weighted sum of its inputs, to assign significance to the inputs with regard to the task the algorithm is trying to learn. It then adds a bias to better represent the part of the output that is actually independent of the input. Finally, it feeds an activation function that determines whether and to what extent the computed value should affect the ultimate outcome. A fully connected network is generally organized into layers, with the neurons of one layer being connected solely to those of the immediately preceding and following layers. The layer that receives the external data is the input layer, the layer that produces the outcome is the output layer, and in between them are zero or more hidden layers.

The design of an efficient neural network requires well-chosen nonlinear activation functions, together with a proper optimization of the weights and biases, to minimize the value of a loss function suitably representing the quality of the network prediction. The network architecture (\eg type of network, depth, width of each layer), the meta-parameters (\ie parameters whose value cannot be estimated from data, \eg optimizer, learning rate, batch size) and the quality/size of the dataset are other key ingredients to an efficient learning, that must be carefully crafted to the intended purpose. For the sake of brevity, the reader is referred to \cite{goodfellow2017} for an extended presentation of this topic.

\subsection{Reinforcement learning}

Reinforcement learning (RL) is a subset of machine learning in which an agent learns to solve decision-making problems by earning rewards through trial-and-error interaction with its environment. 
It is mathematically formulated as a Markov decision process in which the agent observes the current environment state $s_t$, takes an action $a_t$ that prompts
the reward received $r_t$ and the transition to the next state $s_{t+1}$, and repeats until the agent is unable to increase some form of cumulative reward.
In practice, this is framed as an optimization problem of maximizing the discounted reward cumulated over a horizon $T$, defined as:

\begin{equation}
	R(\tau)= \displaystyle \sum_{t=0}^{T} \gamma^t r_{t}\,,
\end{equation}

where $\tau = ( s_0, a_0, s_1, a_1, \dots, s_T, a_T)$ is a trajectory of states and actions, and $\gamma \in [0,1]$ is a discount factor that weights the relative importance of present and future rewards. 

\subsection{Policy gradient RL methods}

In policy gradient methods, the agent behavior is modeled after a stochastic policy $\pi_\theta(s,a)$, \textit{i.e.\@} a probability distribution over actions given states. The expected discounted cumulative reward $J(\theta)$ is maximized by gradient ascent on the policy parameters $\theta$, updating at each iteration by a fixed-size step proportional to the policy gradient, whose expression derived in \cite{sutton2000} in the context of "vanilla" policy gradient reads:

\begin{equation}
\label{eq:policy_gradient}
	\nabla_\theta J(\theta)=\underset{\tau\sim\pi_\theta}{\mathbb{E}} \left[\sum_{t=0}^T \nabla_\theta \log  \pi_\theta (s_t, a_t) R(\tau) \right]\,.
\end{equation}

In a deep reinforcement learning context (deep RL or DRL), the policy is represented by a deep neural network whose weights and biases serve as free parameters to be optimized. To this end, a stochastic gradient algorithm is used to perform network updates from the policy loss:

\begin{equation}
\label{eq:policy_loss}
	L(\theta)= \underset{\tau\sim\pi_\theta}{\mathbb{E}} \left[\sum_{t=0}^T \log \pi_\theta (s_t, a_t) R(\tau) \right]\,,
\end{equation}

whose gradient is equal to $\nabla_\theta J(\theta)$ (since the gradient operator in (\ref{eq:policy_gradient}) acts only on the log-policy term), provided that the transitions $(s_t, a_t, r_t)$ used to perform the update were obtained under policy $\pi_\theta$. The gradient computation is deferred to the back-propagation algorithm \cite{rumelhart1986} with respect to each weight and bias by the chain rule, one layer at the time from the output to the input layer.

Multiple refinements varying in cost, complexity and purpose have been proposed to steer the policy improvement process in the right direction, whether it be by balancing the trade-off between bias and variance (actor-critic \cite{konda2000actor}, generalized advantage estimate \cite{schulman2015gae}) or by preventing the destructively large policy updates that can cause the agent to fall off the cliff and to restart from a poorly performing state with a locally bad policy (trust-region policy optimization, proximal policy optimization \cite{schulman2017}). We shall not elaborate further on this matter, as the present implementation uses a variant of the vanilla policy gradient update. Hence, the interested reader is instead referred to \cite{sutton2018} and references therein.

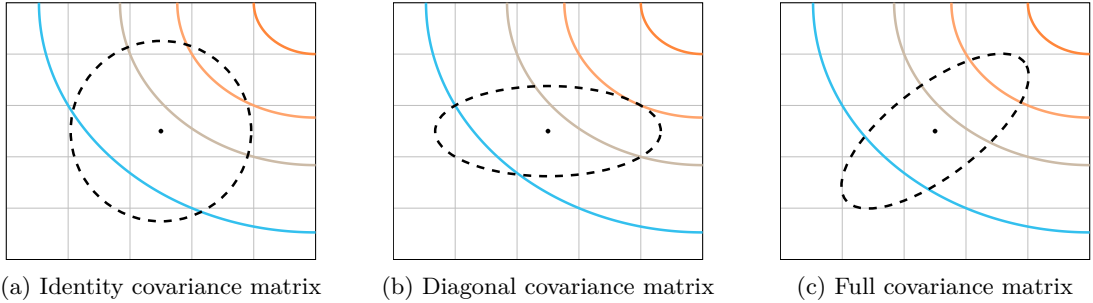
\begin{figure}
\centering
\pgfplotsset{colormap name=custom_map}%
\begin{subfigure}[t]{.3\textwidth}
	\centering
	\begin{tikzpicture}[	scale=0.6, trim axis left, trim axis right,
					global/.style={	circle, fill=black, inner sep=0pt, minimum size=3pt},
					dist/.style={	draw=black, ultra thick, dash pattern=on 5pt}]
		\begin{axis}[	view={0}{90}, clip=false, xmin=-5, xmax=0, ymin=-5, ymax=0, 
					grid=major, ticks=none]

			\addplot [	no markers, raw gnuplot, contour prepared, contour/labels=false,
					point meta min=0, point meta max=25, ultra thick]
				gnuplot {
					set samples 100;
					set isosamples 100;
					set contour base;
					set cntrparam levels discrete 1,5,10,20;
					set style data lines;
					unset surface;
					splot [-5:0] [-5:0] x**2+y**2;
					};
			\node[	global] at (-2.5,-2.5) {};
			\draw[	dist] (-2.5, -2.5) ellipse (2cm and 2cm);
		\end{axis}
	\end{tikzpicture}
    	\caption{Identity covariance matrix}
	\label{fig:identity}
\end{subfigure} \quad
\begin{subfigure}[t]{.3\textwidth}
	\centering
	\begin{tikzpicture}[	scale=0.6, trim axis left, trim axis right,
					global/.style={	circle, fill=black, inner sep=0pt, minimum size=3pt},
					dist/.style={	draw=black, ultra thick, dash pattern=on 5pt}]
		\begin{axis}[	view={0}{90}, clip=false, xmin=-5, xmax=0, ymin=-5, ymax=0, 
					grid=major, ticks=none]

			\addplot [	no markers, raw gnuplot, contour prepared, contour/labels=false,
					point meta min=0, point meta max=25, ultra thick]
				gnuplot {
					set samples 100;
					set isosamples 100;
					set contour base;
					set cntrparam levels discrete 1,5,10,20;
					set style data lines;
					unset surface;
					splot [-5:0] [-5:0] x**2+y**2;
					};
			\node[	global] at (-2.5,-2.5) {};
			\draw[	dist] (-2.5, -2.5) ellipse (2.5cm and 1cm);
		\end{axis}
	\end{tikzpicture}
    	\caption{Diagonal covariance matrix}
	\label{fig:diagonal}
\end{subfigure} \quad
\begin{subfigure}[t]{.3\textwidth}
	\centering
	\begin{tikzpicture}[	scale=0.6, trim axis left, trim axis right,
					global/.style={	circle, fill=black, inner sep=0pt, minimum size=3pt},
					dist/.style={	draw=black, ultra thick, dash pattern=on 5pt}]
		\begin{axis}[	view={0}{90}, clip=false, xmin=-5, xmax=0, ymin=-5, ymax=0, 
					grid=major, ticks=none]

			\addplot [	no markers, raw gnuplot, contour prepared, contour/labels=false,
					point meta min=0, point meta max=25, ultra thick]
				gnuplot {
					set samples 100;
					set isosamples 100;
					set contour base;
					set cntrparam levels discrete 1,5,10,20;
					set style data lines;
					unset surface;
					splot [-5:0] [-5:0] x**2+y**2;
					};
			\node[	global] at (-2.5,-2.5) {};
			\draw[	dist, rotate around={37.5:(-2.5,-2.5)}] (-2.5, -2.5) ellipse (2.5cm and 1cm);
		\end{axis}
	\end{tikzpicture}
    	\caption{Full covariance matrix}
	\label{fig:full}
\end{subfigure}
\caption{\textbf{Iso-density lines for multivariate normal laws} with identity, diagonal and full covariance matrices.}
\label{fig:ellipses}
\end{figure} 

\subsection{Evolution strategies}

Evolution strategies (ES) are another family of stochastic search algorithm 
that can learn an optimal parametrization by emulating organic evolution principles, without knowledge of the performance gradient. At each iteration $g$ (called generation),
the algorithm samples $\lambda$ candidate solutions $(x_1,\dots, x_\lambda)$ from a multivariate normal distribution
$\mathcal{N}(\bm{m}^g,\bm{C}^g)$ 
 with mean $\bm{m}^g$ and covariance matrix $\bm{C}^g$, evaluates the cost function at the candidate solutions, and uses a weighted recombination of the $\mu$ best individuals to update the search distribution for the next generation. 
Simply put, the mean is pulled into the direction of the best performing candidates, while the covariance update aims to align the density contour of the sampling distribution with the contour lines of the objective function and thereby the direction of steepest descent.
The range of possible models corresponds to various degrees of sophistication. For instance,
$(\mu,\lambda)$-ES is a rudimentary algorithm relying on identity covariance matrices, \textit{i.e.}, it assumes all variables to have the same variance and to be uncorrelated, which in turn defines an isotropic region of sampling for the next generation (see figure \ref{fig:identity}). Conversely, the covariance matrix adaptation evolutionary strategy (CMA-ES, considered state-of-the-art in evolutionary computations) uses a full covariance matrix to accelerate convergence toward the optimum by exploiting anisotropy in the steepest descent direction (see figure \ref{fig:full}).
Another key aspect lies in the structure of the CMA-ES covariance matrix update:

\begin{equation}
\label{eq:cmaes}
	\bm{C}^{g+1} \leftarrow ( 1 - c_1 - c_{\mu}) \bm{C}^g + c_{\mu} \bm{C}_{\mu}^{}+ c_1 \bm{C}_1^{}
\end{equation}

where the first term represents a soft update from the current covariance matrix, and $c_1$ and $c_{\mu}$ are learning rates set by well-established heuristics and associated to two types of updates termed \emph{rank-1} and \emph{rank-$\mu$}. The rank-$\mu$ update includes information about the best individuals of the current generation, while the rank-1 update adds correlation information across consecutive generations via 
a so-called evolution path storing the average update direction (in a way such that correlated updates sum up but decorrelated updates cancel each other out). These three contributions combined ultimately allow CMA-ES to fast search from limited populations of individuals at each generation, without compromising the evaluation of the next covariance matrix, as thoroughly described in \cite{hansen2016}.

\section{Policy-based optimization (PBO)}
\label{section:pbo}

We review below the main features of our proposed policy-based optimization (PBO) algorithm, and point out the key conceptual similarities and differences with respect to the methods introduced in section \ref{section:pg-es}. In order to provide common ground between all approaches, we refer from now on to each new set of evaluation as a \emph{generation} $g$, and to each evaluation within a generation as an \emph{individual}. Also, we denote by $n_i$ the number of individuals evaluated at each generation (\textit{i.e.} the number of parallel environments used to collect rewards before performing a network update) and by $d$ the search space dimension (\textit{i.e.} the dimension of the action required by the environment).

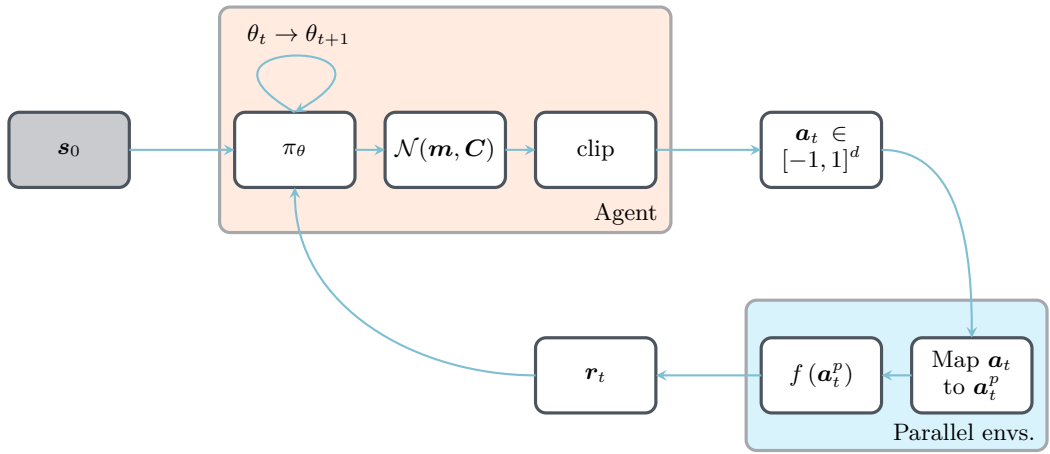
\begin{figure}
\centering
\begin{tikzpicture}[	every loop/.style={	min distance=20mm,looseness=20},
				fontsize/.style={		font=\footnotesize},
				intnode/.style={		black, fontsize, pos=0.5},
				arrow/.style={		thick,color=mybluegray3, rounded corners,thick,-stealth},
				box/.style={		rectangle,rounded corners, draw=mygray1, very thick, 
								text width=1.6cm,minimum height=1cm,text centered,
								inner sep=0pt, outer sep=0pt, fontsize},
				backbox1/.style={	box, opacity=0.5,text width=6cm,minimum height=3cm},
				backbox2/.style={	box, opacity=0.5,text width=4cm,minimum height=2cm}]

	\node[box,fill=mygray4] 			(s0) 		at (0,0) {$\bm{s}_0$};
	
	\node[backbox1,fill=myorange4]	(agent)	at (5,0.4) {};
	\node[xshift=-0.6cm,yshift=0.25cm]	(agt)		at (agent.south east) {\small Agent};
	\node[box,fill=white] 				(nn) 		at (3,0) {$\pi_\theta$};
	\node[box,fill=white] 				(normal) 	at (5,0) {$\mathcal{N}(\bm{m},\bm{C})$};
	\node[box,fill=white] 				(sclip) 	at (7,0) {clip};
	
	\node[box,fill=white] 				(at) 		at (10,0) {$\bm{a}_t \in [-1,1]^{d}$};
	
	\node[backbox2,fill=myblue4]		(env)		at (11,-3) {};
	\node[xshift=-1.1cm,yshift=0.25cm]	(evt)		at (env.south east) {\small Parallel envs.};
	\node[box,fill=white] 				(atp) 		at (12,-3) {Map $\bm{a}_t$ to $\bm{a}^p_t$};
	\node[box,fill=white] 				(f) 		at (10,-3) {$f \left( \bm{a}^p_t \right)$};
	
	\node[box,fill=white] 				(rwd) 	at (7,-3) {$\bm{r}_t$};
	
	\draw[] 		(nn.north) 	edge	[out=150,in=30, loop,arrow] node[intnode, above] {$\,\,\theta_t \to \theta_{t+1}$} (nn.north);
	
	\draw[arrow] 	(s0) 			to 	[out=0,in=180] 	(nn.west);
	\draw[arrow] 	(nn.east)		to 	[out=0,in=180] 	(normal.west);
	\draw[arrow] 	(normal.east)	to 	[out=0,in=180] 	(sclip.west);
	\draw[arrow] 	(sclip.east)	to 	[out=0,in=180] 	(at.west);
	\draw[arrow] 	(at.east)		to 	[out=0,in=90] 	(atp.north);
	\draw[arrow] 	(atp.west)		to 	[out=180,in=0] 	(f.east);
	\draw[arrow] 	(f.west)		to 	[out=180,in=0] 	(rwd.east);
	\draw[arrow] 	(rwd.west)		to 	[out=180,in=-90] (nn.south);
	
\end{tikzpicture}
\caption{\textbf{Action loop for the PBO method.} At each generation, the same input state $\bm{s_0}$ is provided to the agent, that 
draws a set of actions $\bm{a} \in [-1,1]^d$ from the current probability distribution function, with $d$ the problem dimensionality. The actions are distributed to several parallel environments, and mapped to physical ranges $\bm{a^p}$. The parallel environments then evaluate the cost function $f$ at the physical actions, and returns a set of rewards $\bm{r}$ measuring the quality of the actions taken. Once a sufficient amount of state-action-reward triplets has been collected, the network parameters are updated from the policy loss (\ref{eq:pbo_loss}). The process is repeated until convergence.}
\label{fig:pbo}
\end{figure}

\subsection{Single-step deep reinforcement learning}
\label{section:algorithm}

Policy-based optimization (PBO) is a degenerate policy gradient RL algorithm whose premise is that it
is enough to perform single-step episodes if the policy to be learnt is independent of state, \textit{i.e.} $\pi_\theta(s,a)\equiv \pi_\theta(a)$. This is notably the case in optimization and open-loop control problems (the policy in closed-loop control problems conversely depends on states, and thus requires multiple interactions per episode). The line of thought is as follows: where a standard policy gradient algorithm seeks the optimal $\theta^\star$ such that following $\pi_{\theta^\star}$ maximizes the discounted cumulated reward over an episode, PBO seeks the optimal $\theta^\star$ such that $\bm{a^\star} = \pi_{\theta^\star} (\bm{s_0})$ maximizes the instantaneous reward, with $\bm{s_0}$ being some input state (usually a constant vector) consistently fed to the agent for the optimal policy to eventually embody the optimal transformation from $\bm{s_0}$ to $\bm{a^\star}$. The agent initially implements a random policy determined by its initial set of parameters $\theta_0$, after what it gets only one attempt per episode at finding the optimal. This is illustrated in figure \ref{fig:pbo}, showing the agent draw a population of actions from the current policy, and being incentivized to update the policy parameters for the next population of actions to yield larger rewards. A direct consequence is that PBO uses smaller policy networks (compared to usual agent networks found in other DRL contributions), because the agent is not required to learn a complex state-action relation, but only a transformation from a \emph{constant} input state to a given action.

\subsection{Gradient ascent update rule}

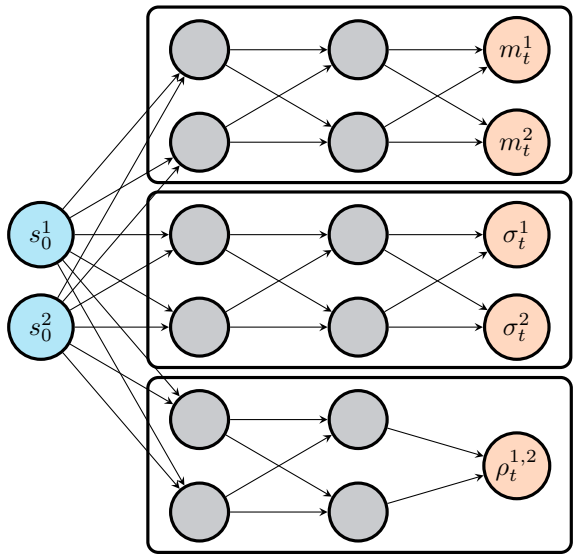
\begin{figure}
\centering
\begin{tikzpicture}[	arrow/.style=		{thick,color=mybluegray1,rounded corners},
				netnode/.style=		{circle, inner sep=0pt, text width=22pt, align=center, very thick},
				inputnode/.style=	{netnode, fill=myblue4, draw=black},
				hiddennode/.style=	{netnode, fill=mygray4, draw=black},
				outputnode/.style=	{netnode, fill=myorange4, 	draw=black},
				signal/.style=		{arrows={-stealth},draw=black}]
	\def\nodedist{35pt}
	\def\layerdist{60pt}
    
	\foreach \y in {1,...,2}
		\node[inputnode] (I\y) at (0,-\y*\nodedist) {$s^\y_0$};  
		
	\foreach \y in {1,...,2}
		\pgfmathsetmacro\z{\y-1}
		\node[hiddennode] (H1\y) at (\layerdist,\z*\nodedist) {};
	\foreach \y in {3,...,4}
		\pgfmathsetmacro\z{\y-2}
		\node[hiddennode] (H1\y) at (\layerdist, -\z*\nodedist) {};
	\foreach \y in {5,...,6}
		\pgfmathsetmacro\z{\y-2}
		\node[hiddennode] (H1\y) at (\layerdist, -\z*\nodedist) {};
		
	\foreach \y in {1,...,2}
		\pgfmathsetmacro\z{\y-1}
		\node[hiddennode] (H2\y) at (2*\layerdist,\z*\nodedist) {};
	\foreach \y in {3,...,4}
		\pgfmathsetmacro\z{\y-2}
		\node[hiddennode] (H2\y) at (2*\layerdist, -\z*\nodedist) {};
	\foreach \y in {5,...,6}
		\pgfmathsetmacro\z{\y-2}
		\node[hiddennode] (H2\y) at (2*\layerdist, -\z*\nodedist) {};
		
	\foreach \y in {1,...,2}
		\pgfmathsetmacro\z{2-\y}
		\node[outputnode] (O\y) at (3*\layerdist, \z*\nodedist) {$m^\y_t$};
	\foreach \y in {3,...,4}
		\pgfmathsetmacro\z{\y-2}
		\pgfmathsetmacro\u{int(\y-2)}
		\node[outputnode] (O\y) at (3*\layerdist, -\z*\nodedist) {$\sigma^\u_t$};
	\foreach \y in {5}
		\pgfmathsetmacro\z{\y-1.5}
		\node[outputnode] (O\y) at (3*\layerdist, -\z*\nodedist) {$\rho^{1,2}_t$};

	\foreach \dest in {1,...,6}
		\foreach \source in {1,...,2}
			\draw[signal] (I\source) -- (H1\dest);
	\foreach \dest in {1,...,2}
		\foreach \source in {1,...,2}
			\draw[signal] (H1\source) -- (H2\dest);
	\foreach \dest in {3,...,4}
		\foreach \source in {3,...,4}
			\draw[signal] (H1\source) -- (H2\dest);
	\foreach \dest in {5,...,6}
		\foreach \source in {5,...,6}
			\draw[signal] (H1\source) -- (H2\dest);
	\foreach \dest in {1,...,2}
		\foreach \source in {1,...,2}
			\draw[signal] (H2\source) edge (O\dest);
	\foreach \dest in {3,...,4}
		\foreach \source in {3,...,4}
			\draw[signal] (H2\source) edge (O\dest);
	\foreach \dest in {5}
		\foreach \source in {5,...,6}
			\draw[signal] (H2\source) edge (O\dest);
			
	\draw[very thick, rounded corners] ($(H11.south west) + (-0.4,-0.25)$) rectangle ($(O1.north east) + (0.4,0.25)$);
	\draw[very thick, rounded corners] ($(H14.south west) + (-0.4,-0.25)$) rectangle ($(O3.north east) + (0.4,0.25)$);
	\draw[very thick, rounded corners] ($(H16.south west) + (-0.4,-0.25)$) rectangle ($(O5.north east) + (0.4,0.85)$);
	
\end{tikzpicture}
\caption{\textbf{Policy networks used in PBO to map states to policy.}  Three separate networks are used for the prediction of mean, standard deviation, and correlation parameters. All activation functions are hyperbolic tangents, except for the output layers of the $\bm{\sigma}$ and $\bm{\rho}$ networks, which uses sigmoid (please refer to section \ref{section:covariance} for additional details). Orthogonal weights initialization is used throughout the networks, with a unit gain for all layers except the output layers, for which the gain is set to \num{1e-2}. In practice, all three networks are trained separately.}
\label{fig:network_pbo}
\end{figure}

In practice, PBO draws actions from a probability density function. Here, we use a $d$-dimensional multivariate normal distribution $\mathcal{N}(\bm{m},\bm{C})$ with mean $\bm{m}$ and full covariance matrix $\bm{C}$. 
As shown in figure \ref{fig:network_pbo}, three independent neural networks are used to output the necessary mean, standard deviation, and correlation information, using hyperbolic tangent and sigmoid activation functions on the output layers to constrain all values in their respective adequate ranges (see section \ref{section:covariance} for more details).
Actions are then drawn in $[-1,1]^d$ by clipping (a series of numerical experiments indicates that soft-limiting transfer functions such as hyperbolic tangent or soft clipping are generally not beneficial and yield, in certain cases, slow convergence and numerical instabilities), before being mapped to their relevant physical ranges $\bm{a^p}$ (a step deferred to the environment as being problem-specific), as illustrated in figure \ref{fig:pbo}. Finally, the Adam algorithm \cite{kingma2014} runs stochastic gradient ascent on the policy parameters using the modified loss:

\begin{equation}
\label{eq:pbo_loss}
	L(\theta)= \underset{a \sim \pi_\theta}{\mathbb{E}} \Big[\log \pi_\theta (a) \, \hat{R} (a) \Big], \text{ with } \hat{R} (a) =  \max \Big(\frac{r(a) - \mu_r}{\sigma_r}, 0 \Big).
 \end{equation}

In the latter expression, $\mu_r$ (resp. $\sigma_r$) is the reward average (resp. standard deviation) over the current generation. The PBO loss is thus formally identical to (\ref{eq:policy_loss}), to the exception of the clipped generation-wise whitened reward substituted for the discounted cumulative reward.
The rationale for this choice is as follows: as is customary in DRL, the discounted cumulative reward is approximated by the advantage function, that measures the improvement (if positive, otherwise the lack thereof) associated with taking action $a$ in state $s$ compared to taking the average over all possible actions. Because a PBO trajectory consists of a single state-action pair (hence (\ref{eq:pbo_loss}) drops the sum over $t$), the discount factor can be set to $\gamma = 1$, in which case the advantage reduces to the reward, as further explained in \cite{ghraieb2020}. The present normalization to zero mean and unit standard deviation introduces bias but reduces variance, and thus the number of actions needed to estimate the expected value. Finally, the max allows discarding negative-advantage actions, that may destabilize learning when performing multiple mini-batch gradient steps using the same data (as each step drives the policy further away from the initial policy).

\subsection{Off-policy updates}
\label{section:off_policy}

Accurately computing the expected value in the policy loss (\ref{eq:pbo_loss}) requires sampling a large number of state-action-reward triplets before the algorithm can proceed to update the agent parameters. At each generation, a set of actions drawn from the current policy $\pi_{\theta}$ is thus distributed to $n_i$ environments running in parallel, each of which computes a reward associated to its input action, and provides it back to the agent. This can repeat until the agent has collected a sufficient number of state-action-reward triplets, Still, in many cases, it is not tractable to use a large value of $n_i$ because computing the reward can be a computationally-intensive task (all the more so when it requires solving high-dimensional discretization of partial differential equation systems), hence the number of state-action-reward triplets available from the current policy is generally limited. Similarly to CMA-ES, PBO therefore improves the reliability of the loss evaluation by incorporating data available from several previous generations. Updating policy $\pi_\theta$ with samples generated under previous policies needs to be accounted for in the loss expression (\ref{eq:pbo_loss}). For samples generated under a policy $\pi_b$, the off-policy loss is written as:

\begin{equation}
\label{eq:pbo_loss_off}
	L_\text{off}(\theta)= \underset{a \sim \pi_b}{\mathbb{E}} \Big[ \frac{\pi_\theta (a)}{\pi_b (a)} \, \hat{R} (a) \Big],
 \end{equation}
 
where $\frac{\pi_\theta (a)}{\pi_b (a)}$ is the \emph{importance} term \cite{degris2013}. The resulting gradient of the objective function is therefore:
 
\begin{equation}
\label{eq:pbo_grad}
	\nabla_\theta J(\theta) = \underset{a \sim \pi_b}{\mathbb{E}} \left[ \frac{\pi_\theta (a)}{\pi_b (a)} \nabla_\theta \log  \pi_\theta (a) \hat{R} (a) \right]\,,
\end{equation}

and the original loss is recovered for $\pi_b = \pi_\theta$. Yet, in practice, it was observed that using (\ref{eq:pbo_loss_off}) led to unstable updates in the final steps of the optimization process, in the vicinity of local or global minima, thus considerably degrading the overall performance of the algorithm. The careful study of this issue is deferred to a future contribution, and, in the meantime, the use of a decay parameter $\eta \in [0,1]$ is introduced to give recent generations more weight by exponentially decreasing the reward from previous generations, in conjunction with loss expression (\ref{eq:pbo_loss}). A rule of thumb for the decay factor used in the remaining of this paper is given by:

\begin{equation}
\label{eq:decay}
	\eta = 1 - \text{e}^{-\alpha d},
\end{equation}

with $\alpha>0$ to retain a longer memory of the previous individuals as the problem dimensionality $d$ increases (very much consistent with the idea that more individuals are then needed to build a coherent covariance matrix). The decrease rate is set empirically to $\alpha=0.35$, hence $\eta=0.5$ for $d=2$, $0.82$ for $\eta=d=5$, and $\eta=0.98$ for $d=10$.

In practice, each of the three neural networks are updated for $n_e$ epochs (the number of full passes of the algorithm over the entire data set) using a learning rate $\lambda_r$ and a history of $n_g$ generations, shuffled and organized in $n_b$ mini-batches (whose size are in multiples of $n_i$, the number of individuals sampled at each generation). An important attribute of PBO is that all three networks can use different meta-parameters and network architectures, which we show in the following can substantially impact the convergence rate.

\subsection{Generating valid covariance matrices from neural network outputs}
\label{section:covariance}

Matrices representing correlations between variables must satisfy four basic properties to bear physical significance: (i) all entries must be in $[-1,1]$ (nothing goes beyond perfect correlation or perfect anticorrelation), (ii) all diagonal entries must be equal to 1 (a variable is always perfectly correlated with itself), (iii) the matrix must be symmetric (correlation between variables $i$ and $j$ is equal to
correlation between $j$ and $i$), and (iv) the matrix must be positive semidefinite (PSD, the variance of a weighted sum of the random variables must be positive). It follows that the above naive approach consisting in having a neural network directly output a set of correlation parameters in adequate range is vowed to fail, as there is no guarantee whatsoever that the so-obtained matrix will be PSD. 
In addition, while it is possible on paper to have the neural network repeatedly output correlation coefficients until a PSD matrix is obtained (which amounts to implementing the classical rejection sampling method), this quickly becomes inefficient as the chances of finding a valid matrix are very low for $d>3$. 

PBO overcomes this issue using hypersphere decomposition, a method rooted in risk management theory, that generates valid correlation matrices from a set of angular coordinates on a hypersphere of unit radius \cite{rebonato2011most,rapisarda2007}. The reader interested in a detailed and comprehensive presentation of the method is referred to \cite{numpacharoen2012}. We shall just mention here that the method parameterizes a lower triangular elementary matrix $\bm{B_d} $ with entry:

\begin{equation}
\label{eq:hypersphere}
	b_{ij} = \left\{
		\begin{array}{ll}
			1 										&\text{for } i=j=1 \\
			\cos \varphi_{ij}								&\text{for } i>1, j=1\\
			\cos \varphi_{ij} \prod_{k=1}^{j-1} \sin \varphi_{ik} 	&\text{for } i>1, j<i\\
			\prod_{k=1}^{j-1} \sin \varphi_{ik} 				&\text{for } i>1, j=i\\			
			0										&\text{for } j>i \\
		\end{array}
		\right.
 \end{equation}

from a set of so-called correlative angles $\bm{\varphi} \in [0,\pi]^D$, with $D = \frac{d(d-1)}{2}$ (hence in same number as the correlation parameters). For instance the matrix for $d=4$ reads:

\begin{equation}
	\bm{B_4} = \left[
		\begin{array}{cccc}
			1 				& 0 							& 0 										& 0\\
			\cos \varphi_{2,1} 	& \sin \varphi_{2,1} 				& 0 										& 0\\
			\cos \varphi_{3,1}	& \cos \varphi_{3,2}  \sin \varphi_{3,1} 	& \sin \varphi_{3,2} \sin \varphi_{3,1} 				& 0\\
			\cos \varphi_{4,1}	& \cos \varphi_{4,2}  \sin \varphi_{4,1} 	& \cos \varphi_{4,3} \sin \varphi_{4,2} \sin \varphi_{4,1} 	& \sin \varphi_{4,3} \sin \varphi_{4,2} \sin \varphi_{4,1}\\
		\end{array}
	\right]
 \end{equation}

The product of this matrix with its transpose is then guaranteed to be a valid correlation matrix, as it is symmetric and PSD by construction, with all entries in $[-1,1]$ (since all $B_{ij}$ are products of cosine and sine functions) and unit diagonal \cite{maree2012correcting}.

The retained procedure to efficiently doctor neural network outputs into valid parameterization of a multivariate normal distribution is thus as follows: the first network outputs the mean $\bm{m}$ in $[-1,1]^d$ using a hyperbolic tangent activation function on the output layer. The second network outputs the standard deviations $\bm{\sigma}$ in $[0,1]^d$ using a sigmoid activation function on the output layer. Finally, the third network outputs a set of coefficients $\bm{\rho}$ in $[0,1]^D$, also using a sigmoid activation function on the output layer. Those are mapped into correlative angles $\bm{\varphi} = \pi \bm{\rho}$ and assembled into the above elementary matrix $\bm{B_d}$, after which the covariance matrix is constructed as:

\begin{equation}
	\bm{C} = \bm{S} \left( \bm{B} \bm{B}^t \right) \bm{S},
\end{equation}

with $\bm{S} = \text{diag} (\bm{\sigma}) $.

\subsection{PBO pseudo-code}
\label{section:pseudo_alg}

To sum up the content of previous sections, a pseudo-code for the PBO method is provided in algorithm \ref{alg:pbo}, in complement of figure \ref{fig:pbo}. The maximal number of generations for the algorithm to run is denoted $n_g^\text{max}$.

\begin{algorithm}
\caption{PBO algorithm}
\label{alg:pbo}
\begin{algorithmic}[1]
\Initialize{$\pi_\theta$, $n_i$ parallel environments}
\For{$g=0, n_g^\text{max} -1$}
	\Sample{$n_i$ actions/individuals $a_i \in \left[-1,1\right]^d$ from $\pi_\theta$}
	\For{$i=0, n_i -1$} \Comment{This loop is executed in parallel}
		\Provide{action $a_i$ to environment $i$}
		\Retrieve{reward $r_i$ from environment $i$} \Comment{End of "single-step" episode}
	\EndFor
	\Compute{clipped normalized reward $\hat{R}_g$} \Comment{Modified reward vector for generation $g$}
	\For{$\sigma$, $\rho$ and $\mu$ networks}
		\For{$e=0, n_e-1$} \Comment{$n_e$ can be specific to each network}
			\Shuffle{data from most recent $n_g$ generations} \Comment{$n_g$ can be specific to each network}
			\For{$b=0, n_b-1$} \Comment{$n_b$ can be specific to each network}
				\Generate{mini-batch $b$ from shuffled data}
				\Update{current network with loss (\ref{eq:pbo_loss})} \Comment{$\lambda$ can be specific to each network}
			\EndFor
		\EndFor
	\EndFor
\EndFor
\end{algorithmic}
\end{algorithm}

\subsection{Connection to evolutionary strategies}
\label{section:connection_es}

While intrinsically a single-step policy-gradient algorithm, several PBO features are reminiscent of the ES and CMA-ES algorithms introduced in section \ref{section:pg-es}. The main analogies are as follows:

\begin{itemize}
	\item Both ES and PBO exploit the successive updates of a probability density function to generate new samples, eventually leading to the convergence to local or a global minimum. In both cases, normal distributions (isotropic for ES, multivariate with full covariance matrix for CMA-ES and PBO) are used, although the concept of PBO is not strictly bound to it;
	\item PBO computes the policy loss (\ref{eq:pbo_loss}) using only the positive-advantage actions. This keeps the policy consistent with the collected experience data, and is reminiscent of the elitist selection of individuals performed in some CMA-ES update rules \cite{hansen2016};
	\item PBO uses history of previous generations to update the network parameters, in the same way CMA-ES uses an evolution path to add information about correlations across consecutive generations;
	\item PBO exponentially decays the advantage history of older generation, which is also a well-known feature of CMA-ES, where scaled covariance matrices from past generations are re-used in future updates and the influence of previous steps decays exponentially in the evolution path \cite{hansen2016}.
\end{itemize}

Ultimately, PBO can be thought as an evolution strategy without a specific update rule, in the sense that CMA-ES relies on tailored analytical heuristics to directly compute the updated probability density function parameters from the previous sample informations, while the update rules of PBO rely on the specific neural network updates routines.

\section{Minimization of analytic functions}
\label{section:analytic}

\subsection{Test cases}
\label{section:cases}

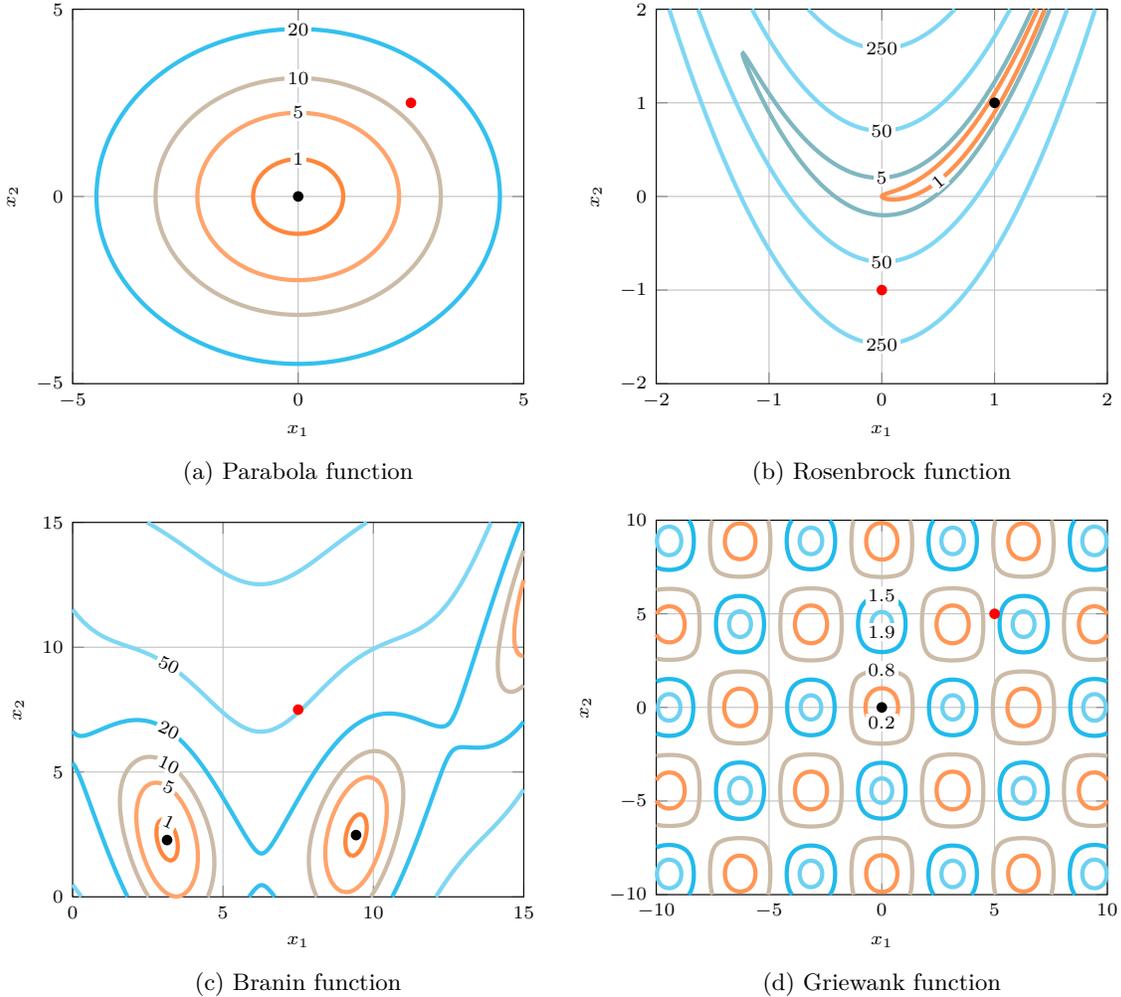
\begin{figure}
\centering
\pgfplotsset{colormap name=custom_map}%
\begin{subfigure}[t]{.4\textwidth}
	\centering
	\begin{tikzpicture}[	trim axis left, trim axis right, font=\scriptsize,
					label/.style={	align=center, fill=white, inner sep=1pt},
					global/.style={	circle, fill=black, inner sep=0pt, minimum size=4pt},
					start/.style={	circle, fill=red, inner sep=0pt, minimum size=4pt}]
		\begin{axis}[	view={0}{90}, clip=false, xmin=-5, xmax=5, ymin=-5, ymax=5, scale=0.875,
					grid=major, xlabel=$x_1$, ylabel=$x_2$, xtick={-5,0,5}, ytick={-5,0,5}]

			\addplot [	no markers, raw gnuplot, contour prepared, contour/labels=false,
					point meta min=0, point meta max=25, ultra thick]
				gnuplot {
					set samples 100;
					set isosamples 100;
					set contour base;
					set cntrparam levels discrete 1,5,10,20;
					set style data lines;
					unset surface;
					splot [-5:5] [-5:5] x**2+y**2;
					};
			\node[	label, rotate=0] at (axis cs: 0,1) {1};
			\node[	label, rotate=0] at (axis cs: 0,2.25) {5};
			\node[	label, rotate=0] at (axis cs: 0,3.15) {10};
			\node[	label, rotate=0] at (axis cs: 0,4.47) {20};
			\node[	global] at (0,0) {};
			\node[	start] at (2.5,2.5) {};
		\end{axis}
	\end{tikzpicture}
    	\caption{Parabola function}
	\label{fig:func_parabola}
\end{subfigure} \qquad \qquad
\begin{subfigure}[t]{.4\textwidth}
	\centering
	\begin{tikzpicture}[	trim axis left, trim axis right, font=\scriptsize,
					label/.style={	align=center, fill=white, inner sep=1pt},
					global/.style={	circle, fill=black, inner sep=0pt, minimum size=4pt},
					start/.style={	circle, fill=red, inner sep=0pt, minimum size=4pt}]
		\begin{axis}[	view={0}{90}, clip=false, xmin=-2, xmax=2, ymin=-2, ymax=2, scale=0.875,
					grid=major, xlabel=$x_1$, ylabel=$x_2$, xtick={-2,-1,0,1,2}, ytick={-2,-1,0,1,2}]

			\addplot [	no markers, raw gnuplot, contour prepared, contour/labels=false,
					point meta min=0, point meta max=10, ultra thick]
				gnuplot {
					set samples 500;
					set isosamples 500;
					set contour base;
					set cntrparam levels discrete 1,5,50,250;
					set style data lines;
					unset surface;
					splot [-2:2] [-2:2] (1-x)**2+100*(y-x**2)**2;
					};
			\node[	label, rotate=40] at (axis cs: 0.5,0.16) {1};
			\node[	label, rotate=0] at (axis cs: 0.0,0.2) {5};
			\node[	label, rotate=0] at (axis cs: 0.0,0.7) {50};
			\node[	label, rotate=0] at (axis cs: 0.0,-0.7) {50};
			\node[	label, rotate=0] at (axis cs: 0.0,1.58) {250};
			\node[	label, rotate=0] at (axis cs: 0.0,-1.58) {250};
			\node[	global] at (1,1) {};
			\node[	start] at (0,-1) {};
		\end{axis}
	\end{tikzpicture}
    	\caption{Rosenbrock function}
	\label{fig:func_rosenbrock}
\end{subfigure}

\medskip

\begin{subfigure}[t]{.4\textwidth}
	\centering
	\begin{tikzpicture}[	trim axis left, trim axis right, font=\scriptsize,
					label/.style={	align=center, fill=white, inner sep=1pt},
					global/.style={	circle, fill=black, inner sep=0pt, minimum size=4pt},
					start/.style={	circle, fill=red, inner sep=0pt, minimum size=4pt}]
		\begin{axis}[	view={0}{90}, clip=false, xmin=0, xmax=15, ymin= 0, ymax=15, scale=0.875,
					grid=major, xlabel=$x_1$, ylabel=$x_2$, xtick={-5,0,5,10}, xtick={0,5,10,15}]

			\addplot [	no markers, raw gnuplot, contour prepared, contour/labels=false,
					point meta min=0, point meta max=25, ultra thick]
				gnuplot {
					set samples 300;
					set isosamples 300;
					set contour base;
					set cntrparam levels discrete 1,5,10,20,50,150;
					set style data lines;
					unset surface;
					splot [0:15] [0:15] (y - 0.129*x**2 + 1.59*x -6)**2 + 9.6*cos(x) + 10;
					};
			\node[	label, rotate=-25] at (axis cs: 3.2,3.0) {1};
			\node[	label, rotate=-25] at (axis cs: 3.2,4.4) {5};
			\node[	label, rotate=-25] at (axis cs: 3.2,5.3) {10};
			\node[	label, rotate=-25] at (axis cs: 3.2,6.7) {20};
			\node[	label, rotate=-25] at (axis cs: 3.2,9.3) {50};
			\node[	global] at (9.42478,2.475) {};
			\node[	global] at (3.14,2.275) {};
			\node[	start] at (7.5,7.5) {};
		\end{axis}
	\end{tikzpicture}
    	\caption{Branin function}
	\label{fig:func_branin}
\end{subfigure} \qquad \qquad
\begin{subfigure}[t]{.4\textwidth}
	\centering
	\begin{tikzpicture}[	trim axis left, trim axis right, font=\scriptsize,
					label/.style={	align=center, fill=white, inner sep=1pt},
					global/.style={	circle, fill=black, inner sep=0pt, minimum size=4pt},
					start/.style={	circle, fill=red, inner sep=0pt, minimum size=4pt}]
		\begin{axis}[	view={0}{90}, clip=false, xmin=-10, xmax=10, ymin=-10, ymax=10, scale=0.875,
					grid=major, xlabel=$x_1$, ylabel=$x_2$]

			\addplot [	no markers, raw gnuplot, contour prepared, contour/labels=false,
					point meta min=0, point meta max=2, ultra thick]
				gnuplot {
					set samples 500;
					set isosamples 500;
					set contour base;
					set cntrparam levels discrete 0.25,0.8,1.5,1.9;
					set style data lines;
					unset surface;
					splot [-10:10] [-10:10] 1+(x**2+y**2)/4000 - cos(x)*cos(y/1.414) ;
					};
			\node[	label, rotate=0] at (axis cs: 0,-0.8) {0.2};
			\node[	label, rotate=0] at (axis cs: 0,2) {0.8};
			\node[	label, rotate=0] at (axis cs: 0,6) {1.5};
			\node[	label, rotate=0] at (axis cs: 0,4) {1.9};
			\node[	global] at (0,0) {};
			\node[	start] at (5,5) {};
 
		\end{axis}
	\end{tikzpicture}
    	\caption{Griewank function}
	\label{fig:func_griewank}
\end{subfigure}
\caption{\textbf{2D analytic functions used as targets for minimization problems.} The global minima and starting points are reported as the black and red dots, respectively.}
\label{fig:functions}
\end{figure} 

This section considers simple minimization problems on a set of analytic functions classically exploited for benchmarking purposes of optimization methods: 
\begin{itemize}
	\item the two-dimensional (2-D) parabola function, whose global minimum is in (0,0), with a search domain equal to $[-5,5]^2$ and a starting point at $(2.5,2.5)$:

	\begin{equation}
	\label{eq:parabola}
		f(x_1,x_2) = x_1^2 + x_2^2\,,
	\end{equation}

	\item the $d$-dimensional ($d$-D) Rosenbrock function, whose global minimum is in $(1,\dots, 1)$ and stands in a very narrow valley notoriously difficult to catch for optimization algorithms (three cases $d=2$, 5 and 10 are tackled for comparison), with a search domain equal to $[-2,2]^d$ and a starting point at $(-1,0)$ in 2-D, and $(0,\dots, 0)$ in 5-D and 10-D:

	\begin{equation}
	\label{eq:rosenbrock}
		f(x_1,\dots x_d) =\sum_{i=1}^{d-1} (1-x_i)^2+100(x_{i+1}-x_i^2)^2\,,
	\end{equation}

	\item the 2-D Branin function, that has two identical global minima at $(\pi,2.275)$ and $(3\pi,2.275)$, with a search domain equal to $[0,15]^2$ and a starting point at $(7.5,7.5)$:

	\begin{equation}
	\label{eq:branin}
		f(x_1,x_2) = \left(x_2 - \frac{5.1}{4 \pi^2} x_1^2 + \frac{5}{\pi} x_1 - 6\right)^2 + 10\left(1 - \frac{1}{8 \pi}\right) \cos(x_1) + 10\,,
	\end{equation}

	\item the 2-D Griewank function, that has multiple widespread, regularly distributed, identical local minima, and only one global minimum at $(0,0)$, with a search domain equal to $[-10,10]^2$ and a starting point at $(5,5)$:

	\begin{equation}
	\label{eq:griewank}
		f(x_1,x_2) = 1 + \frac{x_1^2+x_2^2}{4000} - \cos(x_1)\cos\left(\frac{x_2}{\sqrt{2}}\right)\,.
	\end{equation}
\end{itemize}

The 2-D functions are presented in figure \ref{fig:functions} on their respective domains. In this section, we follow the CMA-ES rules of thumb and set the number of individuals per generation to:

 \begin{equation}
	n_i = \left \lfloor{4 + 3 \ln(d)}\right \rfloor.
\end{equation}

For each case, PBO is benchmarked against our previous single-step PPO-1 algorithm \cite{viquerat2021,ghraieb2020,hachem2020} as well as $(\mu,\lambda)$-ES and CMA-ES algorithms implemented in in-house production codes. To ensure a fair comparison, the initial parameters, number of individuals per generation and starting points are identical for all methods, as indicated in figure \ref{fig:functions}. Moreover, a large initial standard deviation is used by default, to ensure a good exploration of the optimization domain. 

\subsection{Results}
\label{section:resultsbase}

\begin{figure}
\centering
\begin{subfigure}[t]{.45\textwidth}
	\centering
	\begin{tikzpicture}[	trim axis left, trim axis right, font=\scriptsize,
					upper/.style={	name path=upper, smooth, draw=none},
					lower/.style={	name path=lower, smooth, draw=none},]
		\begin{semilogyaxis}[	xmin=0, xmax=50, scale=0.8,
							legend cell align=left, legend pos=south west,
							grid=major, xlabel=generations, ylabel=reward]
			
			\legend{\textsc{ppo-}\oldstylenums{1}, \textsc{es}, \textsc{cmaes}, \textsc{pbo}}
			\addplot [upper, forget plot] 					table[x index=0,y index=3] {data/sphere2d_ppo1_avg.dat}; 
			\addplot [lower, forget plot] 					table[x index=0,y index=2] {data/sphere2d_ppo1_avg.dat}; 
			\addplot [fill=mybluegray4, opacity=0.5, forget plot] 	fill between[of=upper and lower];
			\addplot[draw=mybluegray1, thick, smooth] 		table[x index=0,y index=1] {data/sphere2d_ppo1_avg.dat}; 
			\addplot [upper, forget plot] 					table[x index=0,y index=3] {data/sphere2d_es_avg.dat}; 
			\addplot [lower, forget plot] 					table[x index=0,y index=2] {data/sphere2d_es_avg.dat}; 
			\addplot [fill=mygray4, opacity=0.5, forget plot] 		fill between[of=upper and lower];
			\addplot[draw=mygray1, thick, smooth] 			table[x index=0,y index=1] {data/sphere2d_es_avg.dat}; 
			\addplot [upper, forget plot] 					table[x index=0,y index=3] {data/sphere2d_cmaes_avg.dat}; 
			\addplot [lower, forget plot] 					table[x index=0,y index=2] {data/sphere2d_cmaes_avg.dat}; 
			\addplot [fill=myblue4, opacity=0.5, forget plot] 		fill between[of=upper and lower];
			\addplot[draw=myblue1, thick, smooth] 			table[x index=0,y index=1] {data/sphere2d_cmaes_avg.dat}; 
			\addplot [upper, forget plot] 					table[x index=0,y index=3] {data/sphere2d_pbo_avg.dat}; 
			\addplot [lower, forget plot] 					table[x index=0,y index=2] {data/sphere2d_pbo_avg.dat}; 
			\addplot [fill=myorange4, opacity=0.5, forget plot] 	fill between[of=upper and lower];
			\addplot[draw=myorange1, thick, smooth] 			table[x index=0,y index=1] {data/sphere2d_pbo_avg.dat}; 
			
		\end{semilogyaxis}
	\end{tikzpicture}
    	\caption{2-D parabola function}
	\label{fig:min_sphere2d}
\end{subfigure} \qquad
\begin{subfigure}[t]{.45\textwidth}
	\centering
	\begin{tikzpicture}[	trim axis left, trim axis right, font=\scriptsize,
					upper/.style={	name path=upper, smooth, draw=none},
					lower/.style={	name path=lower, smooth, draw=none},]
		\begin{semilogyaxis}[	xmin=0, xmax=150, scale=0.8, 
							legend cell align=left, legend pos=south west,
							grid=major, xlabel=generations, ylabel=reward]
			
			\legend{\textsc{ppo-}\oldstylenums{1}, \textsc{es}, \textsc{cmaes}, \textsc{pbo}}
			\addplot [upper, forget plot] 					table[x index=0,y index=3] {data/rosenbrock2d_ppo1_avg.dat}; 
			\addplot [lower, forget plot] 					table[x index=0,y index=2] {data/rosenbrock2d_ppo1_avg.dat}; 
			\addplot [fill=mybluegray4, opacity=0.5, forget plot] 	fill between[of=upper and lower];
			\addplot[draw=mybluegray1, thick, smooth] 		table[x index=0,y index=1] {data/rosenbrock2d_ppo1_avg.dat}; 
			\addplot [upper, forget plot] 					table[x index=0,y index=3] {data/rosenbrock2d_es_avg.dat}; 
			\addplot [lower, forget plot] 					table[x index=0,y index=2] {data/rosenbrock2d_es_avg.dat}; 
			\addplot [fill=mygray4, opacity=0.5, forget plot] 		fill between[of=upper and lower];
			\addplot[draw=mygray1, thick, smooth] 			table[x index=0,y index=1] {data/rosenbrock2d_es_avg.dat}; 
			\addplot [upper, forget plot] 					table[x index=0,y index=3] {data/rosenbrock2d_cmaes_avg.dat}; 
			\addplot [lower, forget plot] 					table[x index=0,y index=2] {data/rosenbrock2d_cmaes_avg.dat}; 
			\addplot [fill=myblue4, opacity=0.5, forget plot] 		fill between[of=upper and lower];
			\addplot[draw=myblue1, thick, smooth] 			table[x index=0,y index=1] {data/rosenbrock2d_cmaes_avg.dat}; 
			\addplot [upper, forget plot] 					table[x index=0,y index=3] {data/rosenbrock2d_pbo_avg.dat};
			\addplot [lower, forget plot] 					table[x index=0,y index=2] {data/rosenbrock2d_pbo_avg.dat}; 
			\addplot [fill=myorange4, opacity=0.5, forget plot] 	fill between[of=upper and lower];
			\addplot[draw=myorange1, thick, smooth] 			table[x index=0,y index=1] {data/rosenbrock2d_pbo_avg.dat}; 
			
		\end{semilogyaxis}
	\end{tikzpicture}
    	\caption{2-D Rosenbrock function}
	\label{fig:min_rosenbrock2d}
\end{subfigure}

\bigskip

\begin{subfigure}[t]{.45\textwidth}
	\centering
	\begin{tikzpicture}[	trim axis left, trim axis right, font=\scriptsize,
					upper/.style={	name path=upper, smooth, draw=none},
					lower/.style={	name path=lower, smooth, draw=none},]
		\begin{semilogyaxis}[	xmin=0, xmax=50, scale=0.8, 
							legend cell align=left, legend pos=south west,
							grid=major, xlabel=generations, ylabel=reward]
			
			\legend{\textsc{ppo-}\oldstylenums{1}, \textsc{es}, \textsc{cmaes}, \textsc{pbo}}
			\addplot [upper, forget plot] 					table[x index=0,y index=3] {data/branin_ppo1_avg.dat}; 
			\addplot [lower, forget plot] 					table[x index=0,y index=2] {data/branin_ppo1_avg.dat}; 
			\addplot [fill=mybluegray4, opacity=0.5, forget plot] 	fill between[of=upper and lower];
			\addplot[draw=mybluegray1, thick, smooth] 		table[x index=0,y index=1] {data/branin_ppo1_avg.dat}; 
			\addplot [upper, forget plot] 					table[x index=0,y index=3] {data/branin_es_avg.dat}; 
			\addplot [lower, forget plot] 					table[x index=0,y index=2] {data/branin_es_avg.dat}; 
			\addplot [fill=mygray4, opacity=0.5, forget plot] 		fill between[of=upper and lower];
			\addplot[draw=mygray1, thick, smooth] 			table[x index=0,y index=1] {data/branin_es_avg.dat}; 
			\addplot [upper, forget plot] 					table[x index=0,y index=3] {data/branin_cmaes_avg.dat}; 
			\addplot [lower, forget plot] 					table[x index=0,y index=2] {data/branin_cmaes_avg.dat}; 
			\addplot [fill=myblue4, opacity=0.5, forget plot] 		fill between[of=upper and lower];
			\addplot[draw=myblue1, thick, smooth] 			table[x index=0,y index=1] {data/branin_cmaes_avg.dat}; 
			\addplot [upper, forget plot] 					table[x index=0,y index=3] {data/branin_pbo_avg.dat}; 
			\addplot [lower, forget plot] 					table[x index=0,y index=2] {data/branin_pbo_avg.dat}; 
			\addplot [fill=myorange4, opacity=0.5, forget plot] 	fill between[of=upper and lower];
			\addplot[draw=myorange1, thick, smooth] 			table[x index=0,y index=1] {data/branin_pbo_avg.dat}; 
			
		\end{semilogyaxis}
	\end{tikzpicture}
    	\caption{2-D Branin function}
	\label{fig:min_branin}
\end{subfigure} \qquad
\begin{subfigure}[t]{.45\textwidth}
	\centering
	\begin{tikzpicture}[	trim axis left, trim axis right, font=\scriptsize,
					upper/.style={	name path=upper, smooth, draw=none},
					lower/.style={	name path=lower, smooth, draw=none},]
		\begin{semilogyaxis}[	xmin=0, xmax=50, ymin=0.001, ymax=10.0, scale=0.8, 
							legend cell align=left, legend pos=north east,
							grid=major, xlabel=generations, ylabel=reward, minor tick style={draw=none}]
			
			\legend{\textsc{ppo-}\oldstylenums{1}, \textsc{es}, \textsc{cmaes}, \textsc{pbo}}
			\addplot [upper, forget plot] 					table[x index=0,y index=3] {data/griewank_ppo1_avg.dat}; 
			\addplot [lower, forget plot] 					table[x index=0,y index=2] {data/griewank_ppo1_avg.dat}; 
			\addplot [fill=mybluegray4, opacity=0.5, forget plot] 	fill between[of=upper and lower];
			\addplot[draw=mybluegray1, thick, smooth] 		table[x index=0,y index=1] {data/griewank_ppo1_avg.dat}; 
			\addplot [upper, forget plot] 					table[x index=0,y index=3] {data/griewank_es_avg.dat}; 
			\addplot [lower, forget plot] 					table[x index=0,y index=2] {data/griewank_es_avg.dat}; 
			\addplot [fill=mygray4, opacity=0.5, forget plot] 		fill between[of=upper and lower];
			\addplot[draw=mygray1, thick, smooth] 			table[x index=0,y index=1] {data/griewank_es_avg.dat}; 
			\addplot [upper, forget plot] 					table[x index=0,y index=3] {data/griewank_cmaes_avg.dat}; 
			\addplot [lower, forget plot] 					table[x index=0,y index=2] {data/griewank_cmaes_avg.dat}; 
			\addplot [fill=myblue4, opacity=0.5, forget plot] 		fill between[of=upper and lower];
			\addplot[draw=myblue1, thick, smooth] 			table[x index=0,y index=1] {data/griewank_cmaes_avg.dat}; 
			\addplot [upper, forget plot] 					table[x index=0,y index=3] {data/griewank_pbo_avg.dat}; 
			\addplot [lower, forget plot] 					table[x index=0,y index=2] {data/griewank_pbo_avg.dat}; 
			\addplot [fill=myorange4, opacity=0.5, forget plot] 	fill between[of=upper and lower];
			\addplot[draw=myorange1, thick, smooth] 			table[x index=0,y index=1] {data/griewank_pbo_avg.dat}; 
			
		\end{semilogyaxis}
	\end{tikzpicture}
    	\caption{2-D Griewank function}
	\label{fig:min_griewank}
\end{subfigure}

\bigskip

\begin{subfigure}[t]{.45\textwidth}
	\centering
	\begin{tikzpicture}[	trim axis left, trim axis right, font=\scriptsize,
					upper/.style={	name path=upper, smooth, draw=none},
					lower/.style={	name path=lower, smooth, draw=none},]
		\begin{semilogyaxis}[	xmin=0, xmax=300, scale=0.8, 
							legend cell align=left, legend pos=south west,
							grid=major, xlabel=generations, ylabel=reward]
			
			\legend{\textsc{ppo-}\oldstylenums{1}, \textsc{es}, \textsc{cmaes}, \textsc{pbo}}
			\addplot [upper, forget plot] 					table[x index=0,y index=3] {data/rosenbrock5d_ppo1_avg.dat}; 
			\addplot [lower, forget plot] 					table[x index=0,y index=2] {data/rosenbrock5d_ppo1_avg.dat}; 
			\addplot [fill=mybluegray4, opacity=0.5, forget plot] 	fill between[of=upper and lower];
			\addplot[draw=mybluegray1, thick, smooth] 		table[x index=0,y index=1] {data/rosenbrock5d_ppo1_avg.dat}; 
			\addplot [upper, forget plot] 					table[x index=0,y index=3] {data/rosenbrock5d_es_avg.dat}; 
			\addplot [lower, forget plot] 					table[x index=0,y index=2] {data/rosenbrock5d_es_avg.dat}; 
			\addplot [fill=mygray4, opacity=0.5, forget plot] 		fill between[of=upper and lower];
			\addplot[draw=mygray1, thick, smooth] 			table[x index=0,y index=1] {data/rosenbrock5d_es_avg.dat}; 
			\addplot [upper, forget plot] 					table[x index=0,y index=3] {data/rosenbrock5d_cmaes_avg.dat}; 
			\addplot [lower, forget plot] 					table[x index=0,y index=2] {data/rosenbrock5d_cmaes_avg.dat}; 
			\addplot [fill=myblue4, opacity=0.5, forget plot] 		fill between[of=upper and lower];
			\addplot[draw=myblue1, thick, smooth] 			table[x index=0,y index=1] {data/rosenbrock5d_cmaes_avg.dat}; 
			\addplot [upper, forget plot] 					table[x index=0,y index=3] {data/rosenbrock5d_pbo_avg.dat};
			\addplot [lower, forget plot] 					table[x index=0,y index=2] {data/rosenbrock5d_pbo_avg.dat}; 
			\addplot [fill=myorange4, opacity=0.5, forget plot] 	fill between[of=upper and lower];
			\addplot[draw=myorange1, thick, smooth] 			table[x index=0,y index=1] {data/rosenbrock5d_pbo_avg.dat}; 
			
		\end{semilogyaxis}
	\end{tikzpicture}
    	\caption{5-D Rosenbrock function}
	\label{fig:min_rosenbrock5d}
\end{subfigure}  \qquad
\begin{subfigure}[t]{.45\textwidth}
	\centering
	\begin{tikzpicture}[	trim axis left, trim axis right, font=\scriptsize,
					upper/.style={	name path=upper, smooth, draw=none},
					lower/.style={	name path=lower, smooth, draw=none},]
		\begin{semilogyaxis}[	xmin=0, xmax=600, scale=0.8,
							legend cell align=left, legend pos=south west,
							grid=major, xlabel=generations, ylabel=reward]
			
			\legend{\textsc{ppo-}\oldstylenums{1}, \textsc{es}, \textsc{cmaes}, \textsc{pbo}}
			\addplot [upper, forget plot] 					table[x index=0,y index=3] {data/rosenbrock10d_ppo1_avg.dat}; 
			\addplot [lower, forget plot] 					table[x index=0,y index=2] {data/rosenbrock10d_ppo1_avg.dat}; 
			\addplot [fill=mybluegray4, opacity=0.5, forget plot] 	fill between[of=upper and lower];
			\addplot[draw=mybluegray1, thick, smooth] 		table[x index=0,y index=1] {data/rosenbrock10d_ppo1_avg.dat}; 
			\addplot [upper, forget plot] 					table[x index=0,y index=3] {data/rosenbrock10d_es_avg.dat}; 
			\addplot [lower, forget plot] 					table[x index=0,y index=2] {data/rosenbrock10d_es_avg.dat}; 
			\addplot [fill=mygray4, opacity=0.5, forget plot] 		fill between[of=upper and lower];
			\addplot[draw=mygray1, thick, smooth] 			table[x index=0,y index=1] {data/rosenbrock10d_es_avg.dat}; 
			\addplot [upper, forget plot] 					table[x index=0,y index=3] {data/rosenbrock10d_cmaes_avg.dat}; 
			\addplot [lower, forget plot] 					table[x index=0,y index=2] {data/rosenbrock10d_cmaes_avg.dat}; 
			\addplot [fill=myblue4, opacity=0.5, forget plot] 		fill between[of=upper and lower];
			\addplot[draw=myblue1, thick, smooth] 			table[x index=0,y index=1] {data/rosenbrock10d_cmaes_avg.dat}; 
			\addplot [upper, forget plot] 					table[x index=0,y index=3] {data/rosenbrock10d_pbo_avg.dat}; 
			\addplot [lower, forget plot] 					table[x index=0,y index=2] {data/rosenbrock10d_pbo_avg.dat}; 
			\addplot [fill=myorange4, opacity=0.5, forget plot] 	fill between[of=upper and lower];
			\addplot[draw=myorange1, thick, smooth] 			table[x index=0,y index=1] {data/rosenbrock10d_pbo_avg.dat}; 
			
		\end{semilogyaxis}
	\end{tikzpicture}
    	\caption{10-D Rosenbrock function}
	\label{fig:min_rosenbrock10d}
\end{subfigure}

\caption{\textbf{Minimization problems on analytic functions}, using PBO, PPO-1, ES and CMAES. To ensure a fair comparison, the initial parameters and starting points of the three methods are identical, and the same number of individuals per generation is used for the four methods.}
\label{fig:functions_min}
\end{figure}
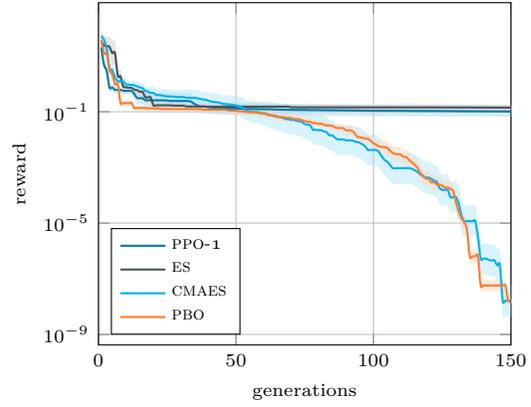
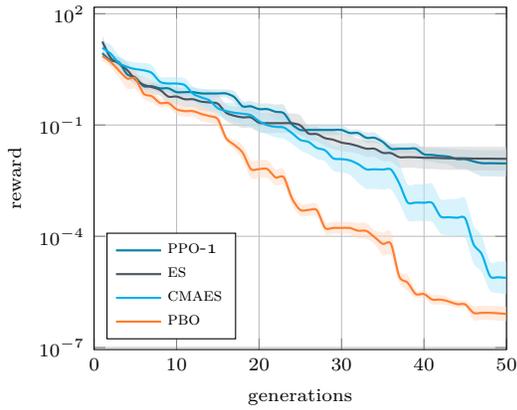
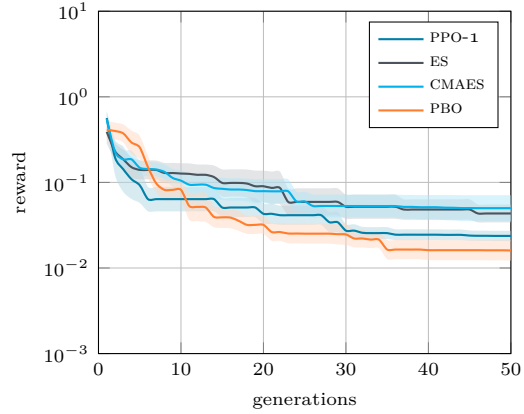
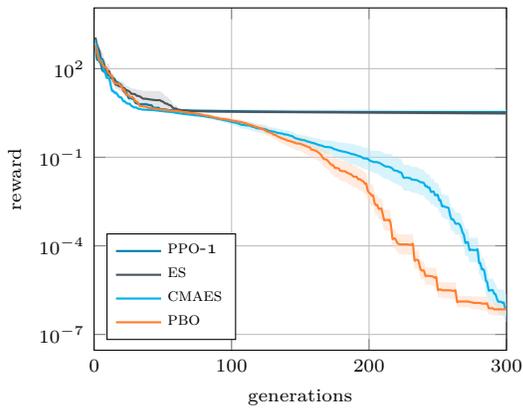
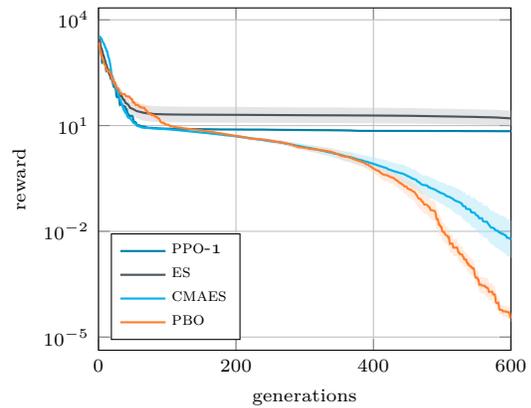

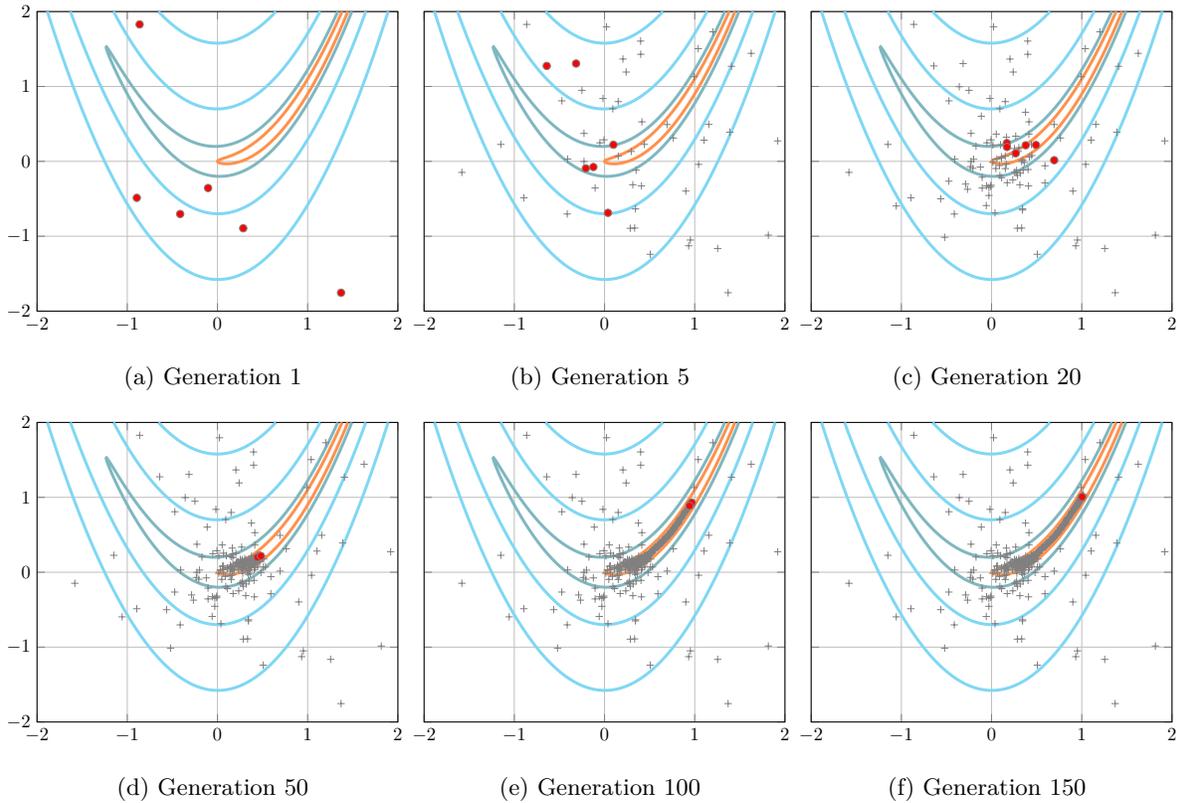
\begin{figure}
\centering
\pgfplotsset{colormap name=custom_map}%
\begin{subfigure}[t]{.3\textwidth}
	\centering
	\begin{tikzpicture}[	scale=0.7, trim axis left, trim axis right,
					label/.style={	align=center, fill=white, inner sep=1pt, font=\tiny}]
		\begin{axis}[	view={0}{90}, clip=false, xmin=-2, xmax=2, ymin=-2, ymax=2, 
					grid=major, xlabel=$$, ylabel=$$, xtick={-2,-1,0,1,2}, ytick={-2,-1,0,1,2}]

			\addplot [	no markers, raw gnuplot, contour prepared, contour/labels=false,
					point meta min=0, point meta max=10, ultra thick]
				gnuplot {
					set samples 500;
					set isosamples 500;
					set contour base;
					set cntrparam levels discrete 1,5,50,250;
					set style data lines;
					unset surface;
					splot [-2:2] [-2:2] (1-x)**2+100*(y-x**2)**2;
					};
			\addplot +[restrict expr to domain={\coordindex}{0:5}, restrict y to domain=-2:2, restrict x to domain=-2:2,only marks, mark=*, mark options={draw=gray, fill=red}] table[x index={3}, y index={4}] {data/pbo_run_rosenbrock};
		\end{axis}
	\end{tikzpicture}
    	\caption{Generation 1}
	\label{fig:pbo_opt_run1}
\end{subfigure} \quad
\begin{subfigure}[t]{.3\textwidth}
	\centering
	\begin{tikzpicture}[	scale=0.7, trim axis left, trim axis right,
					label/.style={	align=center, fill=white, inner sep=1pt, font=\tiny}]
		\begin{axis}[	view={0}{90}, clip=false, xmin=-2, xmax=2, ymin=-2, ymax=2, 
					grid=major, xlabel=$$, ylabel=$$, xtick={-2,-1,0,1,2}, ytick={}, yticklabels={}]

			\addplot [	no markers, raw gnuplot, contour prepared, contour/labels=false,
					point meta min=0, point meta max=10, ultra thick]
				gnuplot {
					set samples 500;
					set isosamples 500;
					set contour base;
					set cntrparam levels discrete 1,5,50,250;
					set style data lines;
					unset surface;
					splot [-2:2] [-2:2] (1-x)**2+100*(y-x**2)**2;
					};
			\addplot +[restrict expr to domain={\coordindex}{0:59}, restrict y to domain=-2:2, restrict x to domain=-2:2,only marks, mark=+, mark options={draw=gray, fill=gray}] table[x index={3}, y index={4}] {data/pbo_run_rosenbrock};
			\addplot +[restrict expr to domain={\coordindex}{60:65}, restrict y to domain=-2:2, restrict x to domain=-2:2,only marks, mark=*, mark options={draw=gray, fill=red}] table[x index={3}, y index={4}] {data/pbo_run_rosenbrock};
		\end{axis}
	\end{tikzpicture}
    	\caption{Generation 5}
	\label{fig:pbo_opt_run5}
\end{subfigure} \quad
\begin{subfigure}[t]{.3\textwidth}
	\centering
	\begin{tikzpicture}[	scale=0.7, trim axis left, trim axis right,
					label/.style={	align=center, fill=white, inner sep=1pt, font=\tiny}]
		\begin{axis}[	view={0}{90}, clip=false, xmin=-2, xmax=2, ymin=-2, ymax=2, 
					grid=major, xlabel=$$, ylabel=$$, xtick={-2,-1,0,1,2}, ytick={}, yticklabels={}]

			\addplot [	no markers, raw gnuplot, contour prepared, contour/labels=false,
					point meta min=0, point meta max=10, ultra thick]
				gnuplot {
					set samples 500;
					set isosamples 500;
					set contour base;
					set cntrparam levels discrete 1,5,50,250;
					set style data lines;
					unset surface;
					splot [-2:2] [-2:2] (1-x)**2+100*(y-x**2)**2;
					};
			\addplot +[restrict expr to domain={\coordindex}{0:119}, restrict y to domain=-2:2, restrict x to domain=-2:2,only marks, mark=+, mark options={draw=gray, fill=gray}] table[x index={3}, y index={4}] {data/pbo_run_rosenbrock};
			\addplot +[restrict expr to domain={\coordindex}{120:125}, restrict y to domain=-2:2, restrict x to domain=-2:2,only marks, mark=*, mark options={draw=gray, fill=red}] table[x index={3}, y index={4}] {data/pbo_run_rosenbrock};
		\end{axis}
	\end{tikzpicture}
    	\caption{Generation 20}
	\label{fig:pbo_opt_run20}
\end{subfigure}

\medskip

\begin{subfigure}[t]{.3\textwidth}
	\centering
	\begin{tikzpicture}[	scale=0.7, trim axis left, trim axis right,
					label/.style={	align=center, fill=white, inner sep=1pt, font=\tiny}]
		\begin{axis}[	view={0}{90}, clip=false, xmin=-2, xmax=2, ymin=-2, ymax=2, 
					grid=major, xlabel=$$, ylabel=$$, xtick={-2,-1,0,1,2}, ytick={-2,-1,0,1,2}]

			\addplot [	no markers, raw gnuplot, contour prepared, contour/labels=false,
					point meta min=0, point meta max=10, ultra thick]
				gnuplot {
					set samples 500;
					set isosamples 500;
					set contour base;
					set cntrparam levels discrete 1,5,50,250;
					set style data lines;
					unset surface;
					splot [-2:2] [-2:2] (1-x)**2+100*(y-x**2)**2;
					};
			\addplot +[restrict expr to domain={\coordindex}{0:299}, restrict y to domain=-2:2, restrict x to domain=-2:2,only marks, mark=+, mark options={draw=gray, fill=red}] table[x index={3}, y index={4}] {data/pbo_run_rosenbrock};
			\addplot +[restrict expr to domain={\coordindex}{300:305}, restrict y to domain=-2:2, restrict x to domain=-2:2,only marks, mark=*, mark options={draw=gray, fill=red}] table[x index={3}, y index={4}] {data/pbo_run_rosenbrock};
		\end{axis}
	\end{tikzpicture}
    	\caption{Generation 50}
	\label{fig:pbo_opt_run50}
\end{subfigure} \quad
\begin{subfigure}[t]{.3\textwidth}
	\centering
	\begin{tikzpicture}[	scale=0.7, trim axis left, trim axis right,
					label/.style={	align=center, fill=white, inner sep=1pt, font=\tiny}]
		\begin{axis}[	view={0}{90}, clip=false, xmin=-2, xmax=2, ymin=-2, ymax=2, 
					grid=major, xlabel=$$, ylabel=$$, xtick={-2,-1,0,1,2}, ytick={}, yticklabels={}]

			\addplot [	no markers, raw gnuplot, contour prepared, contour/labels=false,
					point meta min=0, point meta max=10, ultra thick]
				gnuplot {
					set samples 500;
					set isosamples 500;
					set contour base;
					set cntrparam levels discrete 1,5,50,250;
					set style data lines;
					unset surface;
					splot [-2:2] [-2:2] (1-x)**2+100*(y-x**2)**2;
					};
			\addplot +[restrict expr to domain={\coordindex}{0:599}, restrict y to domain=-2:2, restrict x to domain=-2:2,only marks, mark=+, mark options={draw=gray, fill=gray}] table[x index={3}, y index={4}] {data/pbo_run_rosenbrock};
			\addplot +[restrict expr to domain={\coordindex}{600:605}, restrict y to domain=-2:2, restrict x to domain=-2:2,only marks, mark=*, mark options={draw=gray, fill=red}] table[x index={3}, y index={4}] {data/pbo_run_rosenbrock};
		\end{axis}
	\end{tikzpicture}
    	\caption{Generation 100}
	\label{fig:pbo_opt_run100}
\end{subfigure} \quad
\begin{subfigure}[t]{.3\textwidth}
	\centering
	\begin{tikzpicture}[	scale=0.7, trim axis left, trim axis right,
					label/.style={	align=center, fill=white, inner sep=1pt, font=\tiny}]
		\begin{axis}[	view={0}{90}, clip=false, xmin=-2, xmax=2, ymin=-2, ymax=2, 
					grid=major, xlabel=$$, ylabel=$$, xtick={-2,-1,0,1,2}, ytick={}, yticklabels={}]

			\addplot [	no markers, raw gnuplot, contour prepared, contour/labels=false,
					point meta min=0, point meta max=10, ultra thick]
				gnuplot {
					set samples 500;
					set isosamples 500;
					set contour base;
					set cntrparam levels discrete 1,5,50,250;
					set style data lines;
					unset surface;
					splot [-2:2] [-2:2] (1-x)**2+100*(y-x**2)**2;
					};
			\addplot +[restrict expr to domain={\coordindex}{0:893}, restrict y to domain=-2:2, restrict x to domain=-2:2, only marks, mark=+, mark options={draw=gray, fill=gray}] table[x index={3}, y index={4}] {data/pbo_run_rosenbrock};
			\addplot +[restrict expr to domain={\coordindex}{894:899}, restrict y to domain=-2:2, restrict x to domain=-2:2, only marks, mark=*, mark options={draw=gray, fill=red}] table[x index={3}, y index={4}] {data/pbo_run_rosenbrock};
		\end{axis}
	\end{tikzpicture}
    	\caption{Generation 150}
	\label{fig:pbo_opt_run150}
\end{subfigure}
\caption{\textbf{Successive generations of the PBO algorithm} during a single minimization run of the 2D Rosenbrock function. The red dots indicate the individuals of the current generation, while gray crosses correspond to the individuals of all previous generations.}
\label{fig:pbo_opt_run}
\end{figure}

In order to emphasize flexibility and generalizability, all benchmarks are tackled without fine-tuning of the algorithm, \ie all runs use the same PBO meta-parameters listed in table \ref{table:parameters}, hence the results documented hereafter should be understood as a baseline measure of performance for which there is ample room for improvement.
For each considered case, we present in figure \ref{fig:functions_min} the evolution of the best individual cost during the optimization process of a given algorithm. Performances are averaged over 10 runs, with standard deviations shown as the light shade around. PBO can be seen to perform extremely well on the parabola, and the 2-D Rosenbrock functions, as it significantly outperforms PPO-1 and $(\mu,\lambda)$-ES (both of which perform remarkably similarly) and generally achieves convergence rates and final cost levels similar to CMA-ES. On the 2-D Branin function, the convergence of the PBO algorithm is faster than that of CMA-ES. The anisotropy of the PBO optimization process is further illustrated in figure \ref{fig:pbo_opt_run} for the 2-D Rosenbrock function: starting in $(0,-1)$ with a large initial variance, the algorithm quickly descends toward the entrance of the narrow valley, in the vicinity of $(0,0)$. After a few tens of generation for exploration, the algorithm figures out the shape of the valley entrance, the search distribution starts to elongate, progresses into the valley, before reaching the global minimum within approximately 100 generations.
PBO performs worst on the Griewank function, as the solutions quickly become trapped by one of the local minima due to the inability to set a suitable step size for the local search process (but all methods considered suffer from the same lack of exploration, and ultimately perform almost identically under the same test conditions). In larger dimensions, PBO shows faster convergence and better performance at intermediate stages (here on the 5-D and 10-D Rosenbrock functions). This experiment confirms the capabilities of PBO to efficiently elongate its research area with respect to the local shape of the cost function, and to converge in moderately large research spaces.

\begin{table}
\centering
\begin{tabular}{cccc}
														\toprule
				& $\bm{m}$	& $\bm{\sigma}$	& $\bm{\rho}$	\\\midrule
$\lambda_r$		& \num{5e-3}	& \num{5e-3}		& \num{1e-3}	\\
$n_g$			& 1			& 8				& 16			\\
$n_e$			& 128		& 16				& 16			\\
$n_b$			& 1			& 4				& 8			\\\midrule
Arch.				& $[2,2,2]$	& $[2,2,2]$		& $[2,2,2]$	\\\bottomrule 
\end{tabular}
\bigskip
\caption{\textbf{Detail of the networks achitecture and PBO meta-parameters.} As mentioned in section \ref{section:algorithm}, $\lambda_r$ is the learning rate, $n_g$ is the number of generations used for learning, $n_e$ is the number of epochs, and $n_b$ is the number of mini-batches. For the architecture, only the sizes of the hidden layers are given.}
\label{table:parameters}
\end{table}

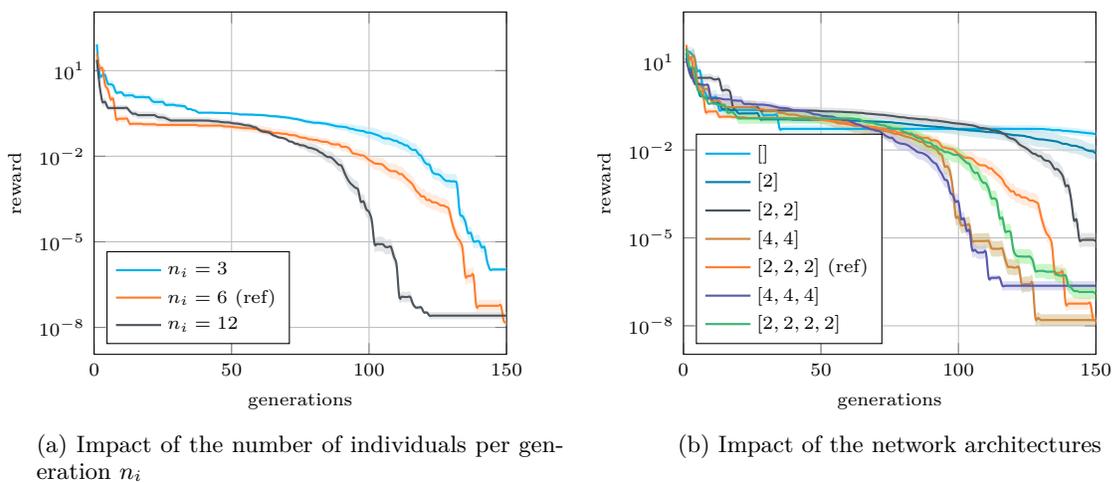
\begin{figure}
\centering
\begin{subfigure}[t]{.45\textwidth}
	\centering
	\begin{tikzpicture}[	trim axis left, trim axis right, font=\scriptsize,
					upper/.style={	name path=upper, smooth, draw=none},
					lower/.style={	name path=lower, smooth, draw=none},]
		\begin{semilogyaxis}[	xmin=0, xmax=150, scale=0.8,
							legend cell align=left, legend pos=south west,
							grid=major, xlabel=generations, ylabel=reward]
			
			\legend{$n_i = 3$, $n_i = 6$ (ref), $n_i = 12$}
			
			\addplot [upper, forget plot] 					table[x index=0,y index=3] {data/rosenbrock2d_pbo_n_3.dat}; 
			\addplot [lower, forget plot] 					table[x index=0,y index=2] {data/rosenbrock2d_pbo_n_3.dat}; 
			\addplot [fill=myblue4, opacity=0.5, forget plot] 	fill between[of=upper and lower];
			\addplot[draw=myblue1, thick, smooth] 			table[x index=0,y index=1] {data/rosenbrock2d_pbo_n_3.dat}; 
			\addplot [upper, forget plot] 					table[x index=0,y index=3] {data/rosenbrock2d_pbo_avg.dat}; 
			\addplot [lower, forget plot] 					table[x index=0,y index=2] {data/rosenbrock2d_pbo_avg.dat}; 
			\addplot [fill=myorange4, opacity=0.5, forget plot] 	fill between[of=upper and lower];
			\addplot[draw=myorange1, thick, smooth] 			table[x index=0,y index=1] {data/rosenbrock2d_pbo_avg.dat}; 
			\addplot [upper, forget plot] 					table[x index=0,y index=3] {data/rosenbrock2d_pbo_n_12.dat}; 
			\addplot [lower, forget plot] 					table[x index=0,y index=2] {data/rosenbrock2d_pbo_n_12.dat}; 
			\addplot [fill=mygray4, opacity=0.5, forget plot] 	fill between[of=upper and lower];
			\addplot[draw=mygray1, thick, smooth] 			table[x index=0,y index=1] {data/rosenbrock2d_pbo_n_12.dat}; 
		
		\end{semilogyaxis}
	\end{tikzpicture}
	\caption{Impact of the number of individuals per generation $n_i$}
	\label{fig:individuals}
\end{subfigure}  \qquad
\begin{subfigure}[t]{.45\textwidth}
	\centering
	\begin{tikzpicture}[	trim axis left, trim axis right, font=\scriptsize,
					upper/.style={	name path=upper, smooth, draw=none},
					lower/.style={	name path=lower, smooth, draw=none},]
		\begin{semilogyaxis}[	xmin=0, xmax=150, scale=0.8,
							legend cell align=left, legend pos=south west,
							grid=major, xlabel=generations, ylabel=reward]
			
			\legend{$[]$\\ $[2]$\\ $[2, 2]$\\ $[4, 4]$\\ $[2, 2, 2]$ (ref)\\ $[4, 4, 4]$\\ $[2, 2, 2, 2]$\\}

			\addplot [upper, forget plot] 					table[x index=0,y index=3] {data/rosenbrock2d_pbo_arch_.dat}; 
			\addplot [lower, forget plot] 					table[x index=0,y index=2] {data/rosenbrock2d_pbo_arch_.dat}; 
			\addplot [fill=myblue4, opacity=0.5, forget plot] 	fill between[of=upper and lower];
			\addplot[draw=myblue1, thick, smooth] 			table[x index=0,y index=1] {data/rosenbrock2d_pbo_arch_.dat}; 
			\addplot [upper, forget plot] 					table[x index=0,y index=3] {data/rosenbrock2d_pbo_arch_2.dat}; 
			\addplot [lower, forget plot] 					table[x index=0,y index=2] {data/rosenbrock2d_pbo_arch_2.dat}; 
			\addplot [fill=mybluegray4, opacity=0.5, forget plot] 	fill between[of=upper and lower];
			\addplot[draw=mybluegray1, thick, smooth] 		table[x index=0,y index=1] {data/rosenbrock2d_pbo_arch_2.dat};
			\addplot [upper, forget plot] 					table[x index=0,y index=3] {data/rosenbrock2d_pbo_arch_2_2.dat}; 
			\addplot [lower, forget plot] 					table[x index=0,y index=2] {data/rosenbrock2d_pbo_arch_2_2.dat}; 
			\addplot [fill=mygray4, opacity=0.5, forget plot] 		fill between[of=upper and lower];
			\addplot[draw=mygray1, thick, smooth] 			table[x index=0,y index=1] {data/rosenbrock2d_pbo_arch_2_2.dat};
			\addplot [upper, forget plot] 					table[x index=0,y index=3] {data/rosenbrock2d_pbo_arch_4_4.dat}; 
			\addplot [lower, forget plot] 					table[x index=0,y index=2] {data/rosenbrock2d_pbo_arch_4_4.dat}; 
			\addplot [fill=mybrown4, opacity=0.5, forget plot] 		fill between[of=upper and lower];
			\addplot[draw=mybrown1, thick, smooth] 			table[x index=0,y index=1] {data/rosenbrock2d_pbo_arch_4_4.dat};
			\addplot [upper, forget plot] 					table[x index=0,y index=3] {data/rosenbrock2d_pbo_avg.dat}; 
			\addplot [lower, forget plot] 					table[x index=0,y index=2] {data/rosenbrock2d_pbo_avg.dat}; 
			\addplot [fill=myorange4, opacity=0.5, forget plot] 	fill between[of=upper and lower];
			\addplot[draw=myorange1, thick, smooth] 			table[x index=0,y index=1] {data/rosenbrock2d_pbo_avg.dat}; 
			\addplot [upper, forget plot] 					table[x index=0,y index=3] {data/rosenbrock2d_pbo_arch_4_4_4.dat}; 
			\addplot [lower, forget plot] 					table[x index=0,y index=2] {data/rosenbrock2d_pbo_arch_4_4_4.dat}; 
			\addplot [fill=mypurple4, opacity=0.5, forget plot] 	fill between[of=upper and lower];
			\addplot[draw=mypurple1, thick, smooth] 			table[x index=0,y index=1] {data/rosenbrock2d_pbo_arch_4_4_4.dat};
			\addplot [upper, forget plot] 					table[x index=0,y index=3] {data/rosenbrock2d_pbo_arch_2_2_2_2.dat}; 
			\addplot [lower, forget plot] 					table[x index=0,y index=2] {data/rosenbrock2d_pbo_arch_2_2_2_2.dat}; 
			\addplot [fill=mygreen4, opacity=0.5, forget plot] 	fill between[of=upper and lower];
			\addplot[draw=mygreen1, thick, smooth] 			table[x index=0,y index=1] {data/rosenbrock2d_pbo_arch_2_2_2_2.dat};
		
		\end{semilogyaxis}
	\end{tikzpicture}
	\caption{Impact of the network architectures}
	\label{fig:architectures}
\end{subfigure}
\caption{\textbf{Sensitivity of the PBO convergence properties to the number of individuals per generation and network architectures.} The reference solutions obtained using the parameters listed in table \ref{table:parameters} and shown in figure \ref{fig:functions_min} are reproduced as the orange curves.}
\label{fig:study}
\end{figure}

We revisit now the 2-D Rosenbrock benchmark and assess the performance sensitivity to the PBO meta-parameters, using the above results (obtained with those meta-parameters listed in table \ref{table:parameters}) as reference.
It can be seen from figure \ref{fig:individuals} that a larger number of individuals per generation $n_i$ leads to a faster convergence. Such a finding is very much consistent with expectations as it proceeds from both a more accurate evaluation of the loss function (\ref{eq:pbo_loss}) and a richer exploration of the search space. Nonetheless, the performance remains surprisingly decent with as little as 3 individuals per generation. The final performance levels seem to saturate around \num{1e-8}, which we believe is a side-effect of the neural network training process. This is a point that deserves further consideration, although the saturation value is small enough that it likely has little to no effect in practical optimization problems.

The architecture of the neural networks also affects performance in a major way, as we show in figure \ref{fig:architectures} that increasing the networks depth and width can make PBO out-perform its reference benchmark (and thus CMA-ES). This means that PBO does indeed exploit the large network parameter state as a proxy to perform efficient optimization, in contrast to just optimizing the bias of the last layer while keeping all weights to zero. Also, PBO being a stochastic method, using deeper networks also substantially increases the performance stability from one run to another. Yet, beyond a certain point, detrimental effects are observed, which can be attributed to vanishing gradients and/or too large parameters states (not shown here). In the same vein, additional numerical experiments (not shown here) conducted on the 5-D and 10-D Rosenbrock functions suggest that these conclusions do not carry over easily to larger dimensional search spaces, and the reference network architecture used in section \ref{section:resultsbase} ends up being a good overall candidate. The general picture to be drawn is that PBO exhibits strong performance and is very promising for use in more applied optimization problems, but that further characterization and fine-tuning are mandatory to outperform more advanced methods on a consistent basis.

\section{Parametric control laws for the Lorenz attractor}
\label{section:lorenz}

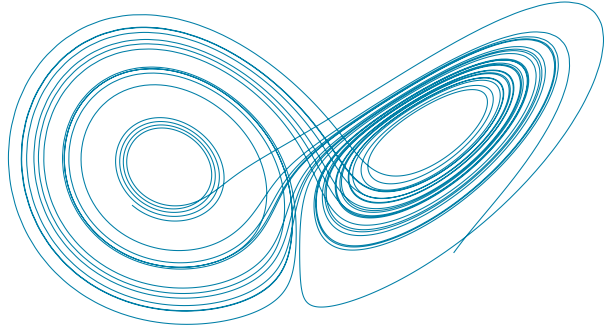
\begin{figure}
\centering
\raisebox{5mm}{
\begin{tikzpicture}
	\begin{axis}[	scale=1.5, hide axis,
				xmin=-25, xmax=25, ymin=-25, ymax=25, zmin=0, zmax=50]
		\addplot3[mark=none, mesh, smooth, color=mybluegray1] table [x index=1, y index=2, z index=3] {data/lorenz_ref.dat}; 
	\end{axis}
\end{tikzpicture}
}
\caption{\textbf{Chaotic sampled solution of the Lorenz attractor}, computed by time-integration of (\ref{eq:lorenz}) with $(\sigma, \rho, \beta) = (10, 28, 8/3)$, from initial conditions $(x_0, y_0, z_0) = (10, 10, 10)$ over 30 time units. The presented view is in the $x-z$ phase plane.}
\label{fig:lorenz}
\end{figure}

This section considers the optimization of parametric control laws for the Lorenz attractor, a simple nonlinear dynamical system representative of thermal convection in a two-dimensional cell \cite{saltzman1962finite}. The set of governing ordinary differential equations reads:

\begin{equation}
\label{eq:lorenz}
\begin{split}
	\dot{x} 	&= \sigma (y - x), \\
	\dot{y}	&= x(\rho - z) - y, \\
	\dot{z}	&= xy - \beta z,
\end{split}
\end{equation}
 
where $\sigma$ is related to the Prandtl number, $\rho$ is a ratio of Rayleigh numbers, and $\beta$ is a geometric factor\footnote{The $\rho$ and $\sigma$ used here are therefore the canonical notations of the Lorenz attractor parameters, and have no link with the standard deviations and correlation parameters used previously in the paper.}. Depending on the values of the triplet $(\sigma, \rho, \beta)$, the solutions to (\ref{eq:lorenz}) may exhibit chaotic behavior, meaning that arbitrarily close initial conditions can lead to significantly different trajectories \cite{lorenz1963}, one common such triplet being $(\sigma, \rho, \beta) = (10, 28, 8/3)$, that leads to the well-known butterfly shape presented in figure \ref{fig:lorenz}. A parametric control law is introduced in the following to alleviate or curb such chaotic behavior, whose design parameters are optimized by PBO with respect each intended control objective.

\subsection{Parametric control law}
\label{section:lorenz_ctrl_laws}

We build here on existing control attempts of the Lorenz system \cite{beintema2020} and add to (\ref{eq:lorenz}) a feedback control on the $y$ variable for the controlled system to be:
 
\begin{equation}
\label{eq:lorenz_ctrl}
\begin{split}
	\dot{x} 	&= \sigma (y - x), \\
	\dot{y}	&= x(\rho - z) - y + u\left(\dot{x}, \dot{y}, \dot{z}\right), \\
	\dot{z}	&= xy - \beta z,
\end{split}
\end{equation}

where $u$ is the feedback velocity defined as:

\begin{equation}
\label{eq:control_law}
	u(\dot{x}, \dot{y}, \dot{z}) = \tanh\left(w_x \dot{x} + w_y \dot{y} + w_z \dot{z} + b \right),
\end{equation}

in a way such that $\left| u \right| < 1$, and $w_x, w_y, w_z$ and $b$ are the true free parameters to optimized (hence $d=4$). The control law (\ref{eq:control_law}) is meant to mimic the output of an artificial neuron, with $w_x, w_y$ and $w_z$ being the weights of the neuron inputs, and $b$ representing its bias. 
We set the initial condition to $(x_0, y_0, z_0) = (10, 10, 10)$, and the attractor parameters to $(\sigma, \rho, \beta) = (10, 28, 8/3)$ for the uncontrolled system to be chaotic. In practice, we use scaled inputs:

\begin{equation}
	\left(\hat{\dot{x}}, \hat{\dot{y}}, \hat{\dot{z}}\right) = \left(\frac{\dot{x}}{x_s}, \frac{\dot{y}}{y_s}, \frac{\dot{z}}{z_s}\right),
\end{equation}

using scaling factors $\left( x_s, y_s, z_s \right) = \left(15,20,40\right)$ representative of the approximate maximal amplitude reached by each variable of the uncontrolled problem, which allows seeking all optimal parameters $w_x^*$, $w_y^*$, $w_z^*$, and $b^*$ in $[-1,1]$ (as required by PBO). In the following, we solve system (\ref{eq:lorenz_ctrl}) using the \verb$odeint$ function of the Scipy package \cite{scipy}. The system is evolved control-free for 5 time units (from $t = -5$ to $t=0$), after which the control kicks in for 25 time units, from $t=0$ to $t=25$. The integration time-step is fixed, and set to $\Delta t=0.01$ time units. All considered cases are tackled with the same reference meta-parameters listed in table \ref{table:parameters} (again to highlight the robustness and versatility of the method before aiming to fine-tune the performance), only the number of individuals per generation $n_i$ is raised to 16 due to the chaotic nature of the system and the limited computational cost required to integrate the problem.

\subsection{Lorenz stabilizer}
\label{section:lorenz_stabilizer}

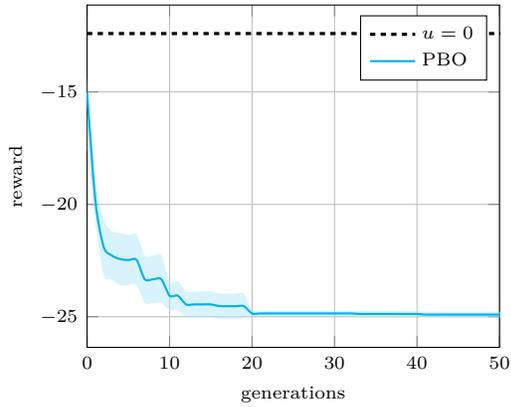
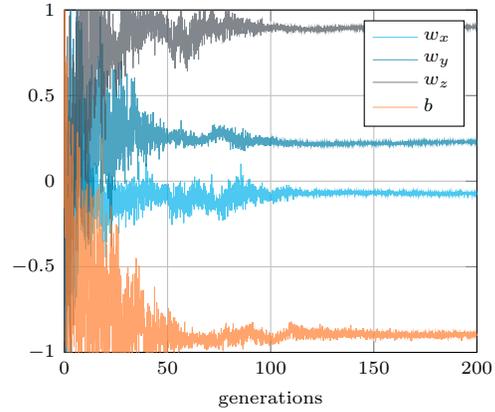
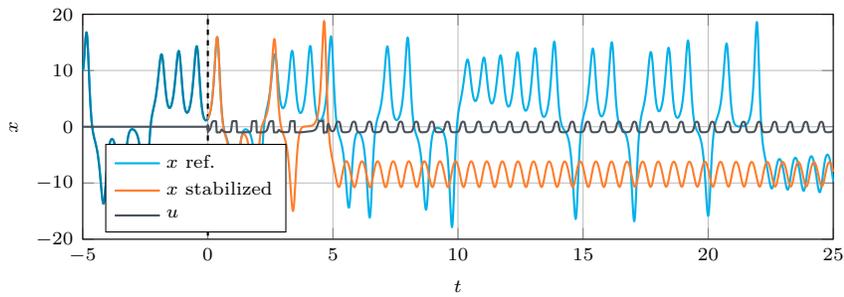
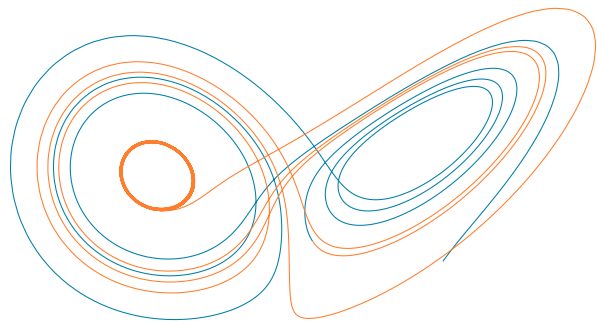
\begin{figure}
\centering
\pgfplotsset{colormap name=custom_map}%
\begin{subfigure}[t]{.45\textwidth}
	\centering
	\begin{tikzpicture}[	trim axis left, trim axis right, font=\scriptsize,
					upper/.style={	name path=upper, smooth, draw=none},
					lower/.style={	name path=lower, smooth, draw=none}]
		\begin{axis}[	xmin=0, xmax=50, scale=0.8, 
					legend cell align=left, legend pos=north east,
					grid=major, xlabel=generations, ylabel=reward]
			
			\legend{$u=0$, PBO}
			
			\addplot[draw=black, very thick, dash pattern=on 2pt, mark=none]	coordinates {(0,-12.4) (200,-12.4)};
			\addplot [upper, forget plot] 					table[x index=0,y index=3] {data/lorenz_stabilizer_avg.dat}; 
			\addplot [lower, forget plot] 					table[x index=0,y index=2] {data/lorenz_stabilizer_avg.dat}; 
			\addplot [fill=myblue4, opacity=0.5, forget plot] 	fill between[of=upper and lower];
			\addplot[draw=myblue1, thick, smooth] 			table[x index=0,y index=1] {data/lorenz_stabilizer_avg.dat }; 
			
		\end{axis}
	\end{tikzpicture}
	\caption{Evolution of the reward function of the Lorenz stabilizer case averaged over 5 runs. The dashed line indicates the control-free reward level}
	\label{fig:lorenz_stabilizer_single}
\end{subfigure}  \qquad
\begin{subfigure}[t]{.45\textwidth}
	\centering
	\begin{tikzpicture}[	trim axis left, trim axis right, font=\scriptsize]
		\begin{axis}[	xmin=0, xmax=200, ymin=-1, ymax=1, scale=0.8, 
					legend cell align=left, legend pos=north east,
					grid=major, xlabel=generations, ylabel={}]
			
			\legend{$w_x$, $w_y$, $w_z$, $b$}
			
			\addplot[draw=myblue1, 		opacity=0.7] table[x expr=\thisrowno{1}/16,y index=3] {data/lorenz_stabilizer_run.dat}; 
			\addplot[draw=mybluegray1, 	opacity=0.7] table[x expr=\thisrowno{1}/16,y index=4] {data/lorenz_stabilizer_run.dat}; 
			\addplot[draw=mygray1, 		opacity=0.7] table[x expr=\thisrowno{1}/16,y index=5] {data/lorenz_stabilizer_run.dat}; 
			\addplot[draw=myorange1, 	opacity=0.7] table[x expr=\thisrowno{1}/16,y index=6] {data/lorenz_stabilizer_run.dat}; 
			
		\end{axis}
	\end{tikzpicture}
	\caption{Evolution of the four control parameters over a single run}
	\label{fig:lorenz_stabilizer_w}
\end{subfigure}

\medskip
\medskip

\begin{subfigure}[t]{.8\textwidth}
	\centering
	\begin{tikzpicture}[	trim axis left, trim axis right, font=\scriptsize]
		\begin{axis}[	xmin=-5, xmax=25, ymin=-20, ymax=20, scale=0.8, width=\textwidth, height=.3\textwidth, scale only axis=true,
					legend cell align=left, legend pos=south west,
					grid=major, xlabel=$t$, ylabel=$x$]
			
			\legend{$x$ ref., $x$ stabilized, $u$}
			
			\addplot[draw=black, 		thick, forget plot, dash pattern=on 2pt] coordinates {(0,-25) (0,25)}; 
			\addplot[draw=myblue1, 		thick, smooth] 			table[x index=0,y index=1] {data/lorenz_ref.dat}; 
			\addplot[draw=mybluegray1, 	thick, smooth, forget plot]	table[x index=0,y index=1] {data/lorenz_stabilizer_plot_init.dat}; 
			\addplot[draw=myorange1, 	thick, smooth]			table[x index=0,y index=1] {data/lorenz_stabilizer_plot_end.dat}; 
			\addplot[draw=mygray1, 		thick, smooth]			table[x index=0,y index=4] {data/lorenz_stabilizer_plot.dat}; 
			
		\end{axis}
	\end{tikzpicture}
	\caption{Evolution of the $x$ component with and without optimal parametric control}
	\label{fig:lorenz_stabilizer_x}
\end{subfigure}

\medskip
\medskip

\begin{subfigure}[t]{.8\textwidth}
	\centering
	\raisebox{5mm}{
	\begin{tikzpicture}
		\begin{axis}[	scale=1.5, hide axis,
					xmin=-25, xmax=25, ymin=-25, ymax=25, zmin=0, zmax=50]
			\addplot3[mark=none, mesh, smooth, color=mybluegray1] 	table [x index=1, y index=2, z index=3] {data/lorenz_stabilizer_plot_init.dat}; 
			\addplot3[mark=none, mesh, smooth, color=myorange1] 		table [x index=1, y index=2, z index=3] {data/lorenz_stabilizer_plot_end.dat }; 
		\end{axis}
	\end{tikzpicture}
	}
	\caption{Plot of the Lorenz attractor with optimized stabilizer control, seen in the $x-z$ plane}
	\label{fig:lorenz_stabilizer_plot}
\end{subfigure}
\caption{\textbf{Results for the Lorenz attractor with optimized stabilizer control}. Given the chaotic nature of the problem, results are provided for a single run only.}
\label{fig:ctrl_stabilizer}
\end{figure}

Small control actuation on the $\dot{y}$ evolution equation is first use to stabilize the Lorenz system in the $x<0$ quadrant (as is done in \cite{beintema2020}) using the reward function:

\begin{equation}
\label{eq:rwd_stabilizer}
	r = \Delta t \sum_{i=0}^{n_t} \left\{x_i < 0\right\},
\end{equation}

where $n_t$ is the total number of time-steps, and $\left\{x_i < 0\right\}=1$ if $x_i < 0$, and 0 otherwise, in a way such that the reward is large if and only if the $x$ coordinate remains within the targeted domain. The reward function is multiplied by $\Delta t$ with the only purpose to make it independent of the time discretization.
For the sake of clarity, the results presented in figure \ref{fig:lorenz_stabilizer_single} pertain to a single run (not an average over runs,) which is because the chaotic behavior of the attractor yields a significantly distorted reward history.
Even though, the PBO algorithm converges after approximately 100 generations (with good reward values are obtained after a few ten generations). The subsequent variations are ascribed to the reward function. Indeed, it is flat by design for any value $x < 0$, and therefore does not promote sharp convergence to a specific value of $x$, although we show in figure \ref{fig:lorenz_stabilizer_w} that all four control parameters converge to well-defined, non-trivial values.
The efficiency of the control is further illustrated in figure \ref{fig:lorenz_stabilizer_x} showing that optimally controlled attractor is successfully confined in the $x<0$ bassin just 5 time units after the control has kicked in.

\subsection{Lorenz oscillator}
\label{section:lorenz_oscillator}

\begin{figure}
\centering
\pgfplotsset{colormap name=custom_map}%
\begin{subfigure}[t]{.45\textwidth}
	\centering
	\begin{tikzpicture}[	trim axis left, trim axis right, font=\scriptsize,
					upper/.style={	name path=upper, smooth, draw=none},
					lower/.style={	name path=lower, smooth, draw=none}]
		\begin{axis}[	xmin=0, xmax=200, scale=0.8, 
					legend cell align=left, legend pos=north east,
					grid=major, xlabel=generations, ylabel=reward]
			
			\legend{$u=0$, PBO}
			
			\addplot[draw=black, very thick, dash pattern=on 2pt, mark=none]	coordinates {(0,-23) (400,-23)};
			\addplot [upper, forget plot] 					table[x index=0,y index=3] {data/lorenz_oscillator_avg.dat }; 
			\addplot [lower, forget plot] 					table[x index=0,y index=2] {data/lorenz_oscillator_avg.dat}; 
			\addplot [fill=myblue4, opacity=0.5, forget plot] 	fill between[of=upper and lower];
			\addplot[draw=myblue1, thick, smooth] 			table[x index=0,y index=1] {data/lorenz_oscillator_avg.dat};
			
		\end{axis}
	\end{tikzpicture}
	\caption{Evolution of the reward function of the Lorenz oscillator case averaged over 5 runs.. The dashed line indicates the control-free reward level}
	\label{fig:lorenz_oscillator_single}
\end{subfigure}  \qquad
\begin{subfigure}[t]{.45\textwidth}
	\centering
	\begin{tikzpicture}[	trim axis left, trim axis right, font=\scriptsize]
		\begin{axis}[	xmin=0, xmax=400, ymin=-1, ymax=1, scale=0.8, 
					legend cell align=left, legend pos=north east,
					grid=major, xlabel=generations, ylabel={}]
			
			\legend{$w_x$, $w_y$, $w_z$, $b$}
			
			\addplot[draw=myblue1, 		opacity=0.7] table[x expr=\thisrowno{1}/16,y index=3] {data/lorenz_oscillator_run.dat }; 
			\addplot[draw=mybluegray1, 	opacity=0.7] table[x expr=\thisrowno{1}/16,y index=4] {data/lorenz_oscillator_run.dat}; 
			\addplot[draw=mygray1, 		opacity=0.7] table[x expr=\thisrowno{1}/16,y index=5] {data/lorenz_oscillator_run.dat}; 
			\addplot[draw=myorange1, 	opacity=0.7] table[x expr=\thisrowno{1}/16,y index=6] {data/lorenz_oscillator_run.dat}; 
			
		\end{axis}
	\end{tikzpicture}
	\caption{Evolution of the four control parameters over a single run}
	\label{fig:lorenz_oscillator_w}
\end{subfigure}

\medskip
\medskip

\begin{subfigure}[t]{.8\textwidth}
	\centering
	\begin{tikzpicture}[	trim axis left, trim axis right, font=\scriptsize]
		\begin{axis}[	xmin=-5, xmax=25, ymin=-20, ymax=20, scale=0.8, width=\textwidth, height=.3\textwidth, scale only axis=true,
					legend cell align=left, legend pos=south west,
					grid=major, xlabel=$t$, ylabel=$x$]
			
			\legend{$x$ ref., $x$ stabilized, $u$}
			
			\addplot[draw=black, 		thick, forget plot, dash pattern=on 2pt] coordinates {(0,-25) (0,25)};
			\addplot[draw=myblue1, 		thick, smooth] 			table[x index=0,y index=1] {data/lorenz_ref.dat};
			\addplot[draw=mybluegray1, 	thick, smooth, forget plot]	table[x index=0,y index=1] {data/lorenz_oscillator_plot_init.dat};
			\addplot[draw=myorange1, 	thick, smooth]			table[x index=0,y index=1] {data/lorenz_oscillator_plot_end.dat};
			\addplot[draw=mygray1, 		thick, smooth]			table[x index=0,y index=7] {data/lorenz_oscillator_plot.dat};
			
		\end{axis}
	\end{tikzpicture}
	\caption{Evolution of the $x$ component with and without optimal parametric control}
	\label{fig:lorenz_oscillator_x}
\end{subfigure}

\medskip
\medskip

\begin{subfigure}[t]{.8\textwidth}
	\centering
	\raisebox{5mm}{
	\begin{tikzpicture}
		\begin{axis}[	scale=1.5, hide axis,
					xmin=-25, xmax=25, ymin=-25, ymax=25, zmin=0, zmax=50]
			\addplot3[mark=none, mesh, smooth, color=mybluegray1] 	table [x index=1, y index=2, z index=3] {data/lorenz_oscillator_plot_init.dat}; 
			\addplot3[mark=none, mesh, smooth, color=myorange1] 		table [x index=1, y index=2, z index=3] {data/lorenz_oscillator_plot_end.dat }; 
		\end{axis}
	\end{tikzpicture}
	}
	\caption{Plot of the Lorenz attractor with optimized oscillator control}
	\label{fig:lorenz_oscillator_plot}
\end{subfigure}
\caption{\textbf{Same as figure \ref{fig:ctrl_stabilizer} for the Lorenz oscillator optimization.}}
\label{fig:ctrl_oscillator}
\end{figure}
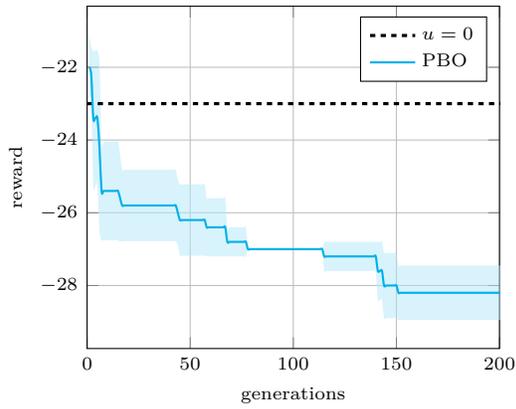
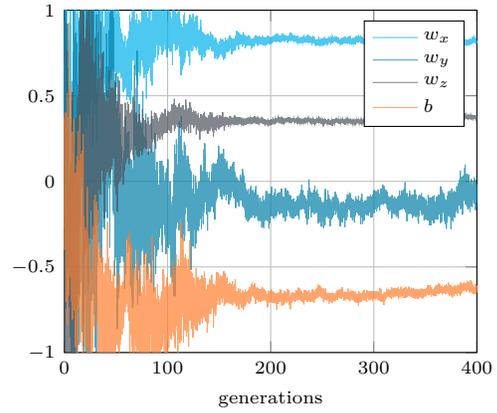
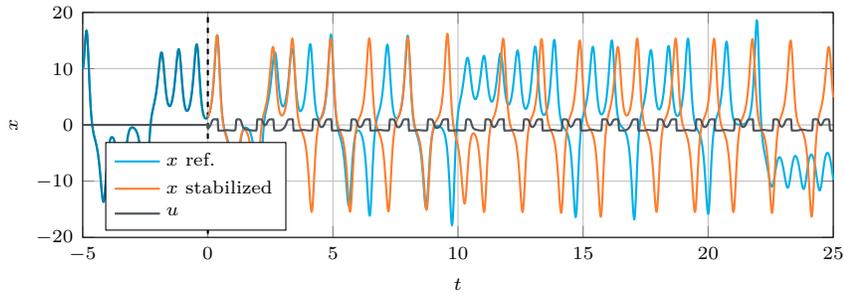
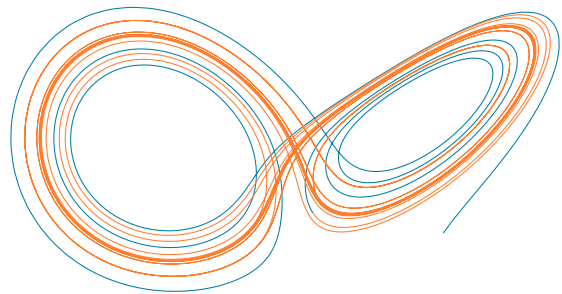

Similar small control actuation on the $\dot{y}$ evolution equation is now used to maximize the number of sign changes of the Lorenz system, as proposed in \cite{beintema2020}. This is done using the following reward function:

\begin{equation}
\label{eq:rwd_oscillator}
	r = \sum_{i=0}^{n_t-1} \left\{ x_i x_{i+1} < 0 \right\},
\end{equation}

in a way such that the reward is large if and only if the $x$ coordinate changes sign in consecutive time steps. This function is much harder to maximize than its stabilizer counterpart (\ref{eq:rwd_stabilizer}) due to its higher sparsity, \textit{i.e.} the larger proportions of actions yielding a zero instantaneous reward.
Such sparsity is the reason why no sharp convergence is found over the course of a single optimization run, as shown in figure \ref{fig:ctrl_oscillator}, but there is a clear diminishing trend and the optimally controlled attractor ultimately exhibits the expected behavior, as it is mostly confined on a narrow orbit that allows it to quickly oscillate between the $x<0$ and the $x>0$ regions.

\section{Conclusion}

This research formally introduces policy-based optimization (PBO), a novel black-box algorithm for optimization and open-loop control problems, at the crossroad of policy gradient methods and evolution strategies. PBO is single-step, meaning that the usual concept of DRL episode is degenerated to a single state-action-reward step. It evolves a multivariate normal search distribution whose parameters (including especially a full covariance matrix) are learnable from neural network outputs.
The method represents significant improvement with respect to our previous single-step PPO-1 algorithm, that samples actions isotropically from a scalar standard deviation (which can be detrimental when the topology of the cost function is distorted), and is shown to outperform classical isotropic ES techniques on the minimization problem of reference analytic functions, up to 10 dimensions. The performance is similar or better to that of CMA-ES, with moderate advantage on the convergence rates obtained in intermediate dimensions, although additional fine-tuning of the method could allow to out-perform its reference benchmark (and thus CMA-ES). PBO is also applied to the optimization of parametric control law for the Lorenz attractor, for which it successfully stabilizes the Lorenz system in a given domain, or conversely enhances the ability of the system to change sign.

Researchers have just begun to gauge the relevance of DRL techniques to assist the design of optimal control strategies. This research weighs in on this issue and shows that PBO holds a high potential as a reliable, go-to black-box
optimizer inheriting from both policy gradients and evolutionary strategy methods.
In this respect, the method can thus benefit from the solid background acquired in evolutionary computations, and from rapid progresses achieved by the DRL community. Despite the present achievements, further development,
characterization and fine-tuning are needed to consolidate the acquired knowledge: multiple refinements can be considered, including extending the scope to deterministic policy gradient techniques \cite{lillicrap2019}, or using importance sampling weights \cite{wang2017} to replace the exponential decay heuristic herein proposed. We certainly welcome such initiatives and purposely make the source code available upon request via a dedicated Github repository \cite{pbo}.

\appendix



\section*{Acknowledgements} 
This work is supported by the Carnot M.I.N.E.S. Institute through the M.I.N.D.S. project.

\bibliographystyle{unsrt}
\bibliography{bib}

\end{document}